\documentclass{article}

\usepackage{amssymb, latexsym}

\usepackage{graphicx}

\usepackage{amsmath}
\usepackage{color}

\numberwithin{equation}{section}

\newtheorem{theorem}{Theorem}[section]

\newtheorem{corollary}[theorem]{Corollary}

\newtheorem{lemma}[theorem]{Lemma}

\newtheorem{proposition}[theorem]{Proposition}

\newtheorem{remark}[theorem]{Remark}

\begin{document}
\title{\large\bf On the  Renyi entropy power  and  the  Gagliardo-Nirenberg-Sobolev inequality on Riemannian manifolds}
\author{\ \ Songzi Li\thanks{Research of S. Li has been supported by NSFC No. 11901569.},\ \ \ \ \  Xiang-Dong Li
\thanks{ Research of X.-D. Li has been supported by NSFC No. 11771430 and Key Laboratory RCSDS, CAS, No. 2008DP173182.} }

\maketitle

\begin{minipage}{120mm}
{\bf Abstract}
In this paper, we prove the concavity of the Renyi entropy power for nonlinear diffusion equation (NLDE)
associated with the Laplacian and the Witten Laplacian on compact
 Riemannian manifolds with non-negative Ricci curvature or $CD(0, m)$-condition and on compact manifolds equipped with time dependent metrics and 
 potentials.  Our results can be regarded as natural extensions of a result due to Savar\'e and Toscani \cite{ST} on the concavity of the Renyi entropy for  NLDE
 on Euclidean spaces. Moreover, we prove that the rigidity models for the Renyi entropy power are the Einstein or quasi-Einstein manifolds and a special $(K, m)$-Ricci flow 
 with Hessian solitons. Inspired by Lu-Ni-Vazquez-Villani \cite{LNVV}, we prove the Aronson-Benilan estimates for  NLDE on compact Riemannian manifolds with 
 $CD(0, m)$-condition.  We also prove the NIW formula which 
 indicates an intrinsic relationship between the second order derivative of the Renyi entropy power $N_p$, the $p$-th Fisher information $I_p$ and 
 the time derivative of the 
 $W$-entropy associated with NLDE.  Finally, we prove the entropy isoperimetric inequality 
 for the Renyi entropy power and the Gagliardo-Nirenberg-Sobolev inequality on complete Riemannian manifolds with non-negative Ricci curvature or $CD(0, m)$-condition and maximal volume growth condition. 
 \end{minipage}


\medskip


\tableofcontents

\section{Introduction}

Let $p>0$. The nonlinear diffusion equation  (NLDE) on $\mathbb{R}^n$ 
\begin{eqnarray}
\partial_t u=\Delta u^p, \label{pme}
\end{eqnarray}
is a nonlinear version of the classical heat equation 
 on $\mathbb{R}^n$
\begin{eqnarray}
\partial_t u=\Delta u. \label{heat1}
\end{eqnarray}
In the  case $p>1$, it is called the {\it 
porous medium equation}, and in the case  $0<p<1$, it is called  the {\it fast diffusion equation}. For $p>1$ and $0<p<1$,  the Cauchy 
problem of $(\ref{pme})$  in weak sense has been well studied in the literature.  Moreover, when $p>1-{2\over n}$, the equation $(\ref{pme})$ preserves 
the mass in the sense that $\int_M u(x, t)dv(x)$ is a constant in $t>0$.  Moreover, when $p\in (1-{2\over n}, 1)$, solutions to  $(\ref{pme})$ are regular
and positive for $t\geq 0$. However, this is not true when $p<1-{2\over n}$, for instance, finite time vanishing may occur. For $p>1$, it is known that solutions to 
$(\ref{pme})$ are at least H\"older continuous. See \cite{DD}.

 In \cite{Ot}, F. Otto proved that the porous medium equation and the fast diffusion equation can be regarded as the 
gradient flow of the Renyi entropy  (see $(\ref{Hp})$ below) on the infinite dimensional $L^2$-Wasserstein space of 
probability measures with finite second moment on Euclidean space equipped with an infinite Riemannian 
metric. See \cite{V1, V2} and reference therein for its extension on Riemannian manifolds. For more background about the porous medium equation and the fast diffusion equation, we refer the readers to 
\cite{Vaz, LNVV, ST}. 

Let  $p\geq 1-{1\over n}$, and let $u$ be a positive and smooth solution to the nonlinear diffusion equation $(\ref{pme})$ with $\int_{\mathbb{R}^n} u(x, t)dx=1$.  Let 
\begin{eqnarray}
H_p(u)={1\over 1-p}\log \int_{\mathbb{R}^n} u^pdx  \label{Hp}
\end{eqnarray}
be the $p$-th R\'enyi entropy associated with the nonlinear diffusion equation $(\ref{pme})$, and define the Renyi entropy power by
\begin{eqnarray*}
N_p(u)=\exp\left(\sigma H_p(u)\right),
\end{eqnarray*}
where $\sigma=p-1+{2\over m}$. In \cite{ST}, Savar\'e and Toscani  proved that  the R\'enyi entropy power $N_p(u(t))$ is concave in $t\in (0, \infty)$, i.e., 
\begin{eqnarray}
{d^2\over dt^2} N_p(u(t)) \leq 0. \label{Np}
\end{eqnarray}
This extends a previous result due to Costa \cite{Cost}, which states that the Shannon entropy power associated with the  heat equation on $\mathbb{R}^n$ is concave. More precisely, let  $u(x, t)$ be  a positive and smooth 
solution to the heat equation $(\ref{heat1})$ on $\mathbb{R}^n$ with $\int_{\mathbb{R}^n} u(x, t)dx=1$,  let 
\begin{eqnarray*}
H(u(t))=-\int_{\mathbb{R}^n} u\log u dx \label{Hu}
\end{eqnarray*}
be the Shannon  entropy associated with $(\ref{heat1})$, and define the Shannon entropy power  by
\begin{eqnarray*}
N(u(t))=e^{{2\over n}H(u(t))}. \label{Nu}
\end{eqnarray*}
Then $N(u(t))$ is concave in $t\in (0, \infty)$, i.e., 
\begin{eqnarray}
{d^2\over dt^2} N(u(t)) \leq 0. \label{N2}
\end{eqnarray}

Using an argument based on the Blachman-Stam inequality \cite{Bla, Stam}, the original proof of $(\ref{N2})$ has been simplified by 
Dembo et al. \cite{Dem1, Dem2}. A direct proof of $(\ref{N2})$ in a strengthened form, with an exact error term, has been obtained by  Villani \cite{V0}, in which  
it was also pointed out that it is possible to extend $(\ref{N2})$ to Riemannian manifolds with non-negative Ricci curvature using the $\Gamma_2$-calculation. 

In our previous work \cite{LL20}, the authors extended Costa's Entropy Power Concavity Inequality
  (EPCI) to the Shannon Entropy Power  for the heat equation associated with the usual Laplacian and the Witten Laplacian 
  on complete Riemannian manifolds with  $CD(K, m)$ curvature-dimension condition and on 
 compact manifolds equipped with $(K, m)$-super Ricci flows. We proved  that the rigidity models for the Shannon entropy power are the Einstein manifolds, quasi-Einstein manifolds or 
 $(K, m)$-Ricci flows equipped with Hessian solitons. Moreover, we proved the NFW formula which indicates an essential relationship between the Shannon entropy power $N$, the $\mathcal{F}$-functional and 
 the time derivative of the 
 $\mathcal{W}$-entropy associated with the conjugate heat equation on Ricci flow, which were introduced  by  G. Perelman \cite{P1}.  As a consequence, we proved that the Shannon entropy power  is convex along the conjugate heat equation on the Ricci flow and 
 the corresponding rigidity models are the shrinking Ricci solitons. This gives us a new understanding of Perelman's entropy formulae on the Ricci  flow \cite{P1}, 
 which plays an important role in the proof of the no local collapsing theorem of Ricci flow for the final resolution of the Poincar\'e conjecture and Thurston's geometrization conjecture. 
 
It is natural to ask the question whether we can extend Savar\'e and Toscani's result to the Renyi entropy power for nonlinear diffusion equations on Riemannian manifolds with suitable curvature 
dimension. The purpose of this paper is to prove the
 Entropy Power Concavity Inequality (EPCI)  of the R\'enyi Entropy Power  for nonlinear diffusion equations
 associated with the usual Laplacian and more general Witten Laplacian on Riemannian manifolds with  $CD(K, m)$ curvature-dimension condition and 
 on manifolds equipped with time dependent metrics and potentials. Indeed,  this part of results have already been proved in 
 our 2017 preprint  \cite{LL17}.  Recently, we have also proved the rigidity theorems and the entropy isoperimetric inequality for the 
 the Renyi entropy power,   and the Gagliardo-Nirenberg-Sobolev inequality as well as the Nash inequality on manifolds with non-negative Ricci curvature or with $CD(0, m)$ curvature-dimension condition and with the maximal volume growth condition.  Moreover, we prove the NIW formula which indicates an essential relationship between the Renyi entropy power $N_p$, the $p$-th Fisher information $I_p$ and 
 the time derivative of the 
 $W$-entropy associated with the nonlinear diffusion equation, which was introduced  by Lu-Ni-Vazquez-Villani \cite{LNVV}. 
 Our work might lead  the readers to recognize  the importance of {\it information-theoretic  approach} in the future study of 
geometric analysis  and Ricci flow. 

\section{Notation and main results}

Let $(M, g)$
be a complete Riemannian manifold, $\phi\in C^2(M)$ and $d\mu=e^{-\phi}dv$, where $v$ is the Riemannian
volume measure on $(M, g)$.
The Witten Laplacian acting on smooth functions is defined by
$$L = \Delta - \nabla\phi\cdot\nabla.$$
For any $u, v \in C_0^\infty(M)$, the integration by parts formula holds
\begin{eqnarray*}
\int_M \langle \nabla u, \nabla v\rangle d\mu=-\int_M L u vd\mu= - \int_M u L v d\mu.
\end{eqnarray*}
Thus, $L$ is the infinitesimal generator of the Dirichlet form
\begin{eqnarray*}
\mathcal{E}(u, v)=\int_M \langle\nabla u, \nabla v\rangle d\mu, \ \ \ \ u, v\in C_0^\infty(M).
\end{eqnarray*}
By It\^o's theory, the Stratonovich SDE on $M$
\begin{eqnarray*}
dX_t =\sqrt{2}U_t\circ dW_t-\nabla\phi(X_t)dt,\ \ \ \ \ \nabla_{\circ dX_t} U_t=0,
\end{eqnarray*}
where $U_t$ is the stochastic parallel transport along the trajectory of $X_t$, with initial data $X_0=x$ and $U_0={\rm Id}_{T_xM}$, defines a diffusion process $X_t$ on $M$ with infinitesimal generator $L$.
Moreover,  the  transition probability density of 
the $L$-diffusion process $X_t$ with respect to $\mu$, i.e., the heat kernel  $p_t(x, y)$   of the Witten Laplacian $L$, is the fundamental solution to the heat equation
\begin{eqnarray}
\partial_t u=Lu.  \label{heq}
\end{eqnarray}

In \cite{BE}, Bakry and Emery proved the generalized Bochner formula 
\begin{eqnarray}
L|\nabla u|^2-2\langle \nabla u, \nabla L u\rangle=2|\nabla^2
u|^2+2Ric(L)(\nabla u, \nabla u), \label{BWF}
\end{eqnarray}
where $u\in C_0^\infty(M)$, $\nabla^2 u$ denotes the Hessian of $u$, $|\nabla^2 u|$ is its Hilbert-Schmidt norm, and 
$$Ric(L)= Ric + \nabla^2\phi$$
is now called the (infinite dimensional) Bakry-Emery Ricci curvature associated with the Witten Laplacian $L$.  For $m\in [n, \infty)$,  the $m$-dimensional Bakry-Emery Ricci curvature associated with the Witten Laplacian $L$  is defined by
$$
Ric_{m, n}(L) = Ric + \nabla^2\phi - {\nabla\phi\otimes\nabla\phi\over m-n}.
$$
In view of this, we have
\begin{eqnarray*}
L|\nabla u|^2-2\langle \nabla u, \nabla L u\rangle \geq {2|Lu|^2\over m}+2Ric_{m, n}(L)(\nabla u, \nabla u).
\end{eqnarray*}
Here we make a convention that $m=n$ if and only if $\phi$ is a constant. By definition, we have
$$Ric(L)=Ric_{\infty, n}(L).$$
Following \cite{BE}, we say that $(M, g, \phi)$ satisfies the curvature-dimension 
$CD(K, m)$-condition for a constant $K\in \mathbb{R}$ and $m\in [n, \infty]$ if and only if
$$Ric_{m, n}(L)\geq Kg.$$
Note that, when $m=n$, $\phi=0$, we have $L=\Delta$ is the usual Laplacian on $(M, g)$, and the $CD(K, n)$-condition holds if and only if the Ricci curvature on $(M, g)$ is bounded from below by $K$, i.e., 
$$Ric\geq Kg.$$

In the case of  Riemannian manifolds with a family of time dependent metrics and potentials, we call $(M, g(t), \phi(t), t\in [0, T])$ a $(K, m)$-super Ricci flow if the metric $g(t)$ and the potential function $\phi(t)$ satisfy 
\begin{eqnarray}
{1\over 2}{\partial g\over \partial t}+Ric_{m, n}(L)\geq Kg, \label{KmsRF}
\end{eqnarray}
where $$L=\Delta_{g(t)}-\nabla_{g(t)}\phi(t)\cdot\nabla_{g(t)}$$ is the time dependent Witten Laplacian on $(M, g(t), \phi(t), t\in [0, T])$, and $K\in \mathbb{R}$ is a constant. When $m=\infty$, i.e., 
if the metric $g(t)$ and the potential function $\phi(t)$ satisfy the following inequality
\begin{eqnarray*}
{1\over 2}{\partial g\over \partial t}+Ric(L)\geq Kg,
\end{eqnarray*}
we call $(M, g(t), \phi(t), t\in [0, T])$   a $(K, \infty)$-super Ricci flow or a $K$-super Perelman Ricci flow. Indeed the $(K, \infty)$-Ricci flow (called also the $K$-Perelman Ricci flow)  
\begin{eqnarray*}
{1\over 2}{\partial g\over \partial t}+Ric(L)=Kg
\end{eqnarray*}
is a natural extension of the modified Ricci flow $\partial_t g=-2Ric(L)$ introduced by Perelman \cite{P1}  as the gradient flow of $\mathcal{F}(g, \phi)=\int_M (R+|\nabla \phi|^2)e^{-\phi}dv$ on $\mathcal{M}\times C^\infty(M)$ 
under the constraint condition that the measure $d\mu=e^{-\phi}dv$ is preserved. For the study of the Li-Yau or Hamilton differential Harnack inequalities, $W$-entropy formulas and related functional inequalities on $(K, m)$ or $(K, \infty)$-super Ricci flows, see \cite{LL15, LL16, LL-AJM,  LL18a, LL18b, LL19a, LL19b}  and references therein. For super Ricci flows on metric and measure spaces, see \cite{Sturm18} and references therein.

We now introduce the Renyi entropy power for the nonlinear diffusion equation  associated with the Witten Laplacian on Riemannian manifolds. Let $u$ be a positive solution to the nonlinear diffusion equation  associated with the Witten Laplacian  on $(M, g, \mu)$
\begin{eqnarray}
\partial_t u=Lu^p, \label{Lup}
\end{eqnarray}
where $0<p<\infty$ and $p\neq 1$.
Following Otto \cite{Ot}, Villani \cite{V1, V2}, Savar\'e and G. Toscani \cite{ST},  for $p\neq 1$,  the $p$-th R\'enyi entropy is defined by
\begin{eqnarray*}
H_p(u)={1\over 1-p}\log \int_M u^pd\mu,
\end{eqnarray*}
and the $p$-Renyi entropy power is defined by 
\begin{eqnarray*}
N_p(u)=\exp\left(\sigma H_p(u)\right),
\end{eqnarray*}
where $\sigma=p-1+{2\over m}$. The $p$-Fisher information is defined by
\begin{eqnarray}
I_p(u):={1\over \int_M u^p d\mu}\int_{\{u>0\}} {|\nabla u^p|^2\over u}d\mu.
\end{eqnarray}
When $p\rightarrow 1$, the Renyi entropy $H_p(u)$, the Renyi entropy power $N_p(u)$ and the $p$-th Fisher information $I_p(u)$ converge to the Shannon entropy $H(u)$, the Shannon entropy power $N(u)$ 
and the Fisher information $I(u)$ associated with the heat equation $(\ref{heq})$. More precisely, 
\begin{eqnarray*}
H(u)&:=&-\int_M u \log u d\mu =\lim\limits_{p\rightarrow 1}H_p(u), \\
N(u)&:=&e^{{2\over m}H(u)}=\lim\limits_{p\rightarrow 1}N_p(u),\\
I(u)&:=&\int_{\{u>0\}} {|\nabla u|^2\over u}d\mu=\lim\limits_{p\rightarrow 1}I_p(u).
\end{eqnarray*}
\medskip

Now we state the main results of this paper.

\begin{theorem}\label{thm1}  Let $M$ be an $n$-dimensional compact Riemannian manifold with $Ric\geq Kg$ for some constant $K\in \mathbb{R}$. Let 
$p\geq 1-{1\over n}$, $\sigma=p-1+{1\over n}$ and  $\kappa=\sigma^{-1}$.
Let $u$ be a positive solution to the nonlinear diffusion equation 
$$\partial_t u=\Delta u^p.$$ 
Then the Renyi entropy power $N_p=N_p(u(t))$ satisfies
\begin{eqnarray}
{d^2 N_p\over dt^2}\leq -{2\sigma K N_p\over \|u\|_p^p} \int_M |\nabla e'(u)|^2 u^pdv. \label{NNp}
\end{eqnarray}
In particular, when $Ric\geq 0$, we have 
\begin{eqnarray}
{d^2N_p\over dt^2}\leq 0. \label{NNp0}
\end{eqnarray}
Moreover,  under the condition $Ric\geq Kg$ (respectively, $Ric\geq 0$), the equality in $(\ref{NNp})$ (respectively,  $(\ref{NNp0})$) holds on $(0, T]$ for some $T>0$ if and only if $(M, g)$ is 
Einstein, i.e., $Ric=Kg$ (respectively, $(M, g)$ is Ricci flat, i.e., $Ric=0$), and $e'(u)$ satisfies 
\begin{eqnarray*}
\Delta e'(u) =\|u\|_p^{-p}\int_M \Delta e'(u) u^pdu=-I_p,\ \ \ \  \nabla^2 e'(u)=-{I_p\over n} g,
\end{eqnarray*} 
where $I_p:={d\over dt}H_{p}(u(t))$ is the $p$-th Fisher information and satisfies the differential equation 
\begin{eqnarray*}
I'_p+\sigma I_p^2+2KI_p=0.
\end{eqnarray*}
When $Ric_{m, n}(L)\geq 0$, it holds
\begin{eqnarray*}
I_p\leq {\kappa \over t}.   \label{Ht2}
\end{eqnarray*}
\end{theorem}

\begin{theorem}\label{thm2} Let $(M, g)$ be an $n$-dimensional compact Riemannian manifold with a potential $\phi\in C^2(M)$ such that $Ric_{m, n}(L)\geq Kg$ for some constant $m\geq n$ and $K\in \mathbb{R}$. Let $p\geq 1-{1\over m}$, $\sigma=p-1+{2\over m}$ and  $\kappa=\sigma^{-1}$.
Let $u$ be a positive solution to the nonlinear diffusion equation 
$$\partial_t u=L u^p.$$ Then the Renyi entropy power  $N_p=N_p(u(t))$  satisfies
\begin{eqnarray}
{d^2N_p\over dt^2}\leq -{2\sigma K N_p\over \|u\|_p^p} \int_M |\nabla e'(u)|^2 u^p d\mu. \label{NNpm}
\end{eqnarray}
In particular, when $Ric\geq 0$, we have 
\begin{eqnarray}
{d^2N_p\over dt^2}\leq 0. \label{NNpm0}
\end{eqnarray}
 Moreover, under the condition $Ric_{m, n}(L)\geq Kg$ (respectively, $Ric_{m, n}(L)\geq 0$), the equality in $(\ref{NNpm})$ (respectively, $(\ref{NNpm0})$) holds on $(0, T]$ for some $T>0$ 
 if and only if $(M, g)$ is quasi-Einstein, i.e., $Ric_{m, n}(L)=Kg$ (respectively, quasi-Ricci flat, i.e., $Ric_{m, n}(L)=0$), and $e'(u)$ satisfies 
\begin{eqnarray}
L e'(u) =\int_M Le'(u) d\gamma=-I_p,\ \ \ \  \nabla^2 e'(u)=-{I_p\over m} g, \ \ \ {m\over m-n}\nabla\phi\cdot\nabla e'(u)=I_p,  \label{rrrr1}
\end{eqnarray}
where $I_p:={d\over dt}H_{p}(u(t))$ is the $p$-th Fisher information and satisfies the differential equation 
\begin{eqnarray*}
I'_p+\sigma I_p^2+2KI_p=0.
\end{eqnarray*}
When $Ric_{m, n}(L)\geq 0$, it holds
\begin{eqnarray}
I_p\leq {\kappa \over t}.   \label{Ht2}
\end{eqnarray}
\end{theorem}

\begin{theorem}\label{thm3} Let $(M, g(t), \phi(t), t\in [0, T])$ be an $n$-dimensional compact Riemannian manifold with potentials $\phi(t)\in C^2(M)$ such that 
\begin{eqnarray}
Ric_{m, n}(L)\geq K_1g,\  \ \ \ \partial_t g\geq 2K_2g, \ \ \ \partial_t \phi={1\over 2 }{\rm Tr} \left({\partial g\over \partial t}\right),  \label{supRF-1}
\end{eqnarray}
where $m\geq n$ and $K, K_1, K_2\in \mathbb{R}$ are constants. Let $p\geq 1-{1\over m}$ and $\sigma=p-1+{2\over m}$. 
Let $u$ be a positive solution to the nonlinear diffusion equation 
$$\partial_t u=L u^p.$$ Then the Renyi entropy power  $N_p=N_p(u(t))$  satisfies
\begin{eqnarray}
{d^2N_p\over dt^2}\leq-{2\sigma K_1 N_p\over \|u\|_p^p} \int_M |\nabla e'(u)|^2 u^p d\mu-{2\sigma K_2 N_p\over \|u\|_p^p}
\int_M |\nabla e'(u)|^2 ud\mu.  \label{NNpmt}
\end{eqnarray}
Moreover, the equality in $(\ref{NNpmt})$ holds on $(0, T]$ for some $T>0$ if and only if 
\begin{eqnarray*}
 Ric_{m, n}(L)=K_1g, \ \ \  \partial_t g=2K_2g,\ \ \ \  \partial_t \phi=nK_2,
\end{eqnarray*}
and $e'(u)$ satisfies    
\begin{eqnarray*}
L e'(u) =\int_M Le'(u) d\gamma=-I_p,\ \ \ \  \nabla^2 e'(u)=-{I_p\over m} g, \ \ \ {m\over m-n}\nabla\phi\cdot\nabla e'(u)=I_p.
\end{eqnarray*}
where $I_p$ is the $p$-th Fisher information and satisfies the differential equation 
\begin{eqnarray*}
I'_p+\sigma I_p^2+2KI_p=0.
\end{eqnarray*}
\end{theorem}

The condition$(\ref{supRF-1})$  in Theorem \ref{thm3}  implies that  $(M, g(t), \phi(t))$ is a $(K, m)$-super Ricci flow with  $K=K_1+K_2$. Conversely, if $(M, g(t), \phi(t))$ is a $(K, m)$-super Ricci flow and with $K_1g\leq Ric_{m, n}(L)\leq 
K_3g$, then $\partial_t g\geq 2K_2 g$ with $K_2=K-K_3$.

The above results extend the result due to Savar\'e and Toscani on the concavity of the Renyi entropy for the nonlinear diffusion equation on Euclidean space to Riemannian manifolds and super Ricci flows.

We will give two proofs for the entropy power concavity inequalities in Theorems \ref{thm1}, \ref{thm2} and \ref{thm3}. 
The first proof uses the similar  idea of  Savar\'e and Toscani \cite{ST} for the proof on the concavity of the Renyi entropy power 
along the nonlinear diffusion equation $(\ref{pme})$ on Euclidean spaces, while the second proof is similar to the one in our previous paper \cite{LL20} for the linear heat equation $(\ref{heq})$ on Riemannian manifolds and super Ricci flows, 
and is based on an explicit formula 
for the second order derivative of the Renyi entropy power  associated with the nonlinear diffusion equation of the usual Laplacian or the Witten Laplacian 
on Riemannian manifolds. In the case $p=1$, such an explicit  formula for the second order derivative of the Shannon entropy power 
associated with the heat equation of the usual Laplacian or the Witten Laplacian  on Riemannian manifolds 
was first proved in our previous paper \cite{LL20}.  We would like point out that each of these two proofs has its own advantage: the first one can be extended to the setting of porous medium equation and 
fast diffusion equation 
on the so-called RCD$(K, N)$ metric measure spaces, which we we will  develop in a forthcoming paper, and the second one enables  us  to prove the rigidity theorems  for the Renyi 
entropy power on compact Riemannin manifolds with $CD(K, m)$-condition or on compact $(K, m)$-super Ricci flows.

The rest of this paper is organized as follows. In Section 2, we prove the dissipation formulae for the Renyi entropy on compact Riehmannian manifolds 
with time dependent or time independent metrics and potentials. In Section 3 we give the first proof of the concavity inequality for the Renyi entropy power in Theorem \ref{thm1}, 
Theorem \ref{thm2} and Theorem \ref{thm3}.  In Section 4, we prove an explicit formula for the second order derivative of the Renyi entropy power on manifolds with the $CD(K, m)$-condition or time dependent metrics and potentials. We will see that 
there is a significant difference for the linear heat equation and the nonlinear diffusion equation on manifolds with time dependent metrics and potentials.  
 In Section 5, we prove an explicit formula 
for the second order derivative of the Renyi entropy power  associated with the nonlinear diffusion equation of the usual Laplacian or the Witten Laplacian 
on Riemannian manifolds and super Ricci flows. In Section 6, we give the second proof of the concavity inequality and prove the rigidity theorems for the Renyi entropy power associated with
  the nonlinear diffusion 
 equation on Riemannian manifolds with the $CD(K, m)$-condition or $(K, m)$-super Ricci flows. In Section 7, we prove the Aronson-Benilan estimates for the nonlinear diffusion equation on compact Riemannian manifolds with $CD(0, m)$-condition. Moreover, we prove the NIW formula which 
gives  an intrinsic relationship between the second order derivative of the Renyi entropy power $N_p$, the $p$-th Fisher information $I_p$ and 
 the time derivative of the 
 $W$-entropy associated with the nonlinear diffusion equation, which was introduced by Lu-Ni-Vazquez-Villani \cite{LNVV}. 
  In Section 8,  we prove an entropy isoperimetric inequality 
 for the Renyi entropy power,  the Gagliardo-Nirenberg-Sobolev inequality and the Nash inequality 
  on complete Riemannian manifolds with non-negative Ricci curvature  or $CD(0, m)$-condition and the maximal volume growth condition.

\medskip

\section{Entropy dissipation formulae for nonlinear diffusion equation}

In this section, we prove two entropy dissipation formulae for the Renyi entropy associated with the nonlinear diffusion equation on 
manifolds with fix metric and potential or time dependent metrics  and potentials.

For  $p>1-2/m$ and $p\neq 1$, let  $$e(r)={1\over p-1}r^p.$$
Define 
\begin{eqnarray*}
E(u):=\int_M e(u(x))d\mu(x),
\end{eqnarray*}
and
\begin{eqnarray*}
E'(u):&=&\int_M |\nabla e'(u)|^2 ud\mu,\\
E''(u):&=&2\int_M \left(|\nabla^2e'(u)|^2+Ric(L)(\nabla e'(u), \nabla e'(u))\right) u^p d\mu\\
& &+2(p-1)\int_M (Le'(u))^2 u^pd\mu+\int_M {\partial g\over \partial t}(\nabla e'(u), \nabla e'(u)) ud\mu.
\end{eqnarray*}

\begin{theorem} \label{thm4} Let $(M, g(t), \phi(t), t\in [0, T])$ be a compact manifold equipped with a family of time dependent metrics  and potentials $(g(t), \phi(t), t\in [0, T])$ which satisfies the conjugate heat equation
\begin{eqnarray}
{\partial \phi\over \partial t}={1\over 2}{\rm Tr}\left(
{\partial g\over \partial t}\right). \label{conjheateq}
\end{eqnarray}
Let $u$ be a positive and smooth solution to the nonlinear diffusion equation $(\ref{Lup})$. Then
\begin{eqnarray*}
{d\over dt}E(u)&=&-E'(u),\\
{d^2\over dt^2}E(u)&=&E''(u).
\end{eqnarray*}
In particular, when $(M, g, \phi)$ is a compact Riemannian manifold with time independent metric and potential, the above formulae holds 
with ${\partial g\over dt}=0$ in $E''(u)$. 
\end{theorem}
{\it Proof}. Let $f(u)=u^p$. Then $f(u)=ue'(u)-e(u)$, $f'(u)=ue''(u)=pu^{p-1}$, and the nonlinear diffusion equation $(\ref{Lup})$ can be rewritten as 
\begin{eqnarray*}
\partial_t u=\nabla_\mu^* (u\nabla e'(u)),
\end{eqnarray*}
where $-\nabla_\mu^*$ denotes the $L^2$-adjoint of $\nabla$ with respect to $\mu$.

Under the conjugate heat equation $(\ref{conjheateq})$, we have
\begin{eqnarray}
{\partial \over \partial t}d\mu=0. \label{invmu}
\end{eqnarray}
Taking time derivative and integrating by parts, we have
\begin{eqnarray*}
{d\over dt}E(u)&=&\int_M e'(u)\partial_t u d\mu\\
&=&\int_M e'(u)\nabla_\mu^*(u\nabla e'(u)) d\mu\\
&=&-\int_M |\nabla e'(u)|^2 ud\mu\\
&=&-E'(u).
\end{eqnarray*}
Taking time derivative again and using the fact (cf. \cite{Lot1, LL15})
\begin{eqnarray*}
\partial_t |\nabla e'(u)|^2=-{\partial g\over \partial t}(\nabla e'(u), \nabla e'(u))+2\langle \nabla e'(u), \nabla \partial_t e'(u)),
\end{eqnarray*}
the  derivative of $E'(u)$  is given by

\begin{eqnarray*}
{d\over dt}E'(u)=\int_M \left[|\nabla e'(u)|^2\partial_t u+2\nabla e'(u)\cdot\nabla (e''(u)\partial_t u)u-{\partial g\over \partial t}(\nabla e'(u), \nabla e'(u))u \right] d\mu.
\end{eqnarray*}
Note that
\begin{eqnarray*}
e''(u)\partial_t u&=&e''(u)\nabla_\mu^*(u\nabla e'(u))\\
&=&ue''(u)Le'(u)+e''(u)\nabla u \cdot\nabla e'(u)\\
&=& f'(u)Le'(u)+|\nabla e'(u)|^2.
\end{eqnarray*}This yields
\begin{eqnarray*}
{d\over dt}E'(u)&=&\int_M |\nabla e'(u)|^2 Lu^pd\mu+2\int_M u\nabla e'(u)\cdot \nabla (f'(u)Le'(u))d\mu\\
& &+2\int_M u\nabla e'(u)\cdot \nabla |\nabla e'(u)|^2 d\mu-\int_M {\partial g\over \partial t}(\nabla e'(u), \nabla e'(u))u d\mu\\
&=&I_1+I_2+I_3+I_4.
\end{eqnarray*}
Integration by parts yields 
\begin{eqnarray*}
I_1=\int_M L|\nabla e'(u)|^2 u^p d\mu=\int_M L|\nabla e'(u)|^2 f(u) d\mu.
\end{eqnarray*}
Splitting $I_2$ into two terms
\begin{eqnarray*}
I_2=2\int_M  u\nabla e'(u)\cdot \nabla (f'(u)Le'(u)) d\mu=I_{21}+I_{22},
\end{eqnarray*}
where
\begin{eqnarray*}
 I_{21} &=&2\int_M u f'(u) \nabla e'(u)\cdot \nabla Le'(u) d\mu,\\
 I_{22} &=& 2\int_M u \nabla e'(u) \cdot \nabla f'(u) Le'(u)d\mu\\
 &=&2\int_M \nabla f(u)\cdot \nabla f'(u)  Le'(u)d\mu.
\end{eqnarray*}
Moreover, we have
\begin{eqnarray*}
I_3&=&2\int_M u\nabla e'(u)\cdot \nabla |\nabla e'(u)|^2 d\mu\\
&=&2\int_M \nabla f(u)\cdot\nabla |\nabla e'(u)|^2d\mu\\
&=&-2\int_M f(u)L|\nabla e'(u)|^2d\mu,
\end{eqnarray*}
and 
\begin{eqnarray*}
I_4=-\int_M {\partial g\over \partial t}(\nabla e'(u), \nabla e'(u))u d\mu.
\end{eqnarray*}
Combining the above calculation together, we can derive that
\begin{eqnarray*}
{d\over dt}E'(u)
&=&-\int_M f(u)[L|\nabla e'(u)|^2-2\nabla e'(u)\cdot\nabla L e'(u)]d\mu\\
& &+2\int_M (u f'(u)-f(u)) \nabla e'(u)\cdot \nabla Le'(u) d\mu\\
& &+2\int_M \nabla f(u)\cdot \nabla f'(u)  Le'(u)d\mu-\int_M {\partial g\over \partial t}(\nabla e'(u), \nabla e'(u))u d\mu.
\end{eqnarray*}
Denote
\begin{eqnarray}
2\Gamma_2(\nabla e'(u), \nabla e'(u)):=L|\nabla e'(u)|^2-2\nabla e'(u)\cdot\nabla L e'(u).
\end{eqnarray}
Integrating by parts, and using 
\begin{eqnarray*}
(rf'(r)-f(r))'=rf''(r), \ \ \ \ \ rf''(r)e''(r)=f'(r)f''(r),
\end{eqnarray*}
we have
\begin{eqnarray*}
{d\over dt}E'(u)&=&-2\int_M f(u)\Gamma_2(\nabla e'(u), \nabla e'(u)))d\mu-2\int_M (u f'(u)-f(u)) |Le'(u)|^2 d\mu\\
& &-2\int_M (uf'(u)-f(u))'\nabla u \cdot \nabla e'(u) Le'(u) d\mu\\
& &+2\int_M \nabla f(u)\cdot \nabla f'(u)  Le'(u)d\mu-\int_M {\partial g\over \partial t}(\nabla e'(u), \nabla e'(u))u d\mu\\
&=&-2\int_M \left[\Gamma_2(\nabla e'(u), \nabla e'(u))+(p-1)(Le'(u))^2 \right]u^pd\mu-\int_M {\partial g\over \partial t}(\nabla e'(u), \nabla e'(u))u d\mu,
\end{eqnarray*}
where in the last step we have used the fact 
\begin{eqnarray*}
& & f(u)=u^p,\ \ \ \ \ u f'(u)-f(u)=(p-1)f,\\
& &uf''(u)\nabla u\cdot\nabla e'(u)=p^2(p-1)u^{2p-3}|\nabla u|^2,\\
& &\nabla f(u)\cdot \nabla f'(u)=f'(u)f''(u)|\nabla u|^2=p^2(p-1)u^{2p-3}|\nabla u|^2.
\end{eqnarray*}
By  the Bakry-Emery-Bochner formula (\cite{BE}), we have
\begin{eqnarray*}
L|\nabla e'(u)|^2-2\nabla e'(u)\cdot\nabla L e'(u)
=2|\nabla^2 e'(u)|^2+2Ric(L)(\nabla e'(u), \nabla e'(u)).
\end{eqnarray*}
 This completes the proof of Theorem \ref{thm4}. \hfill $\square$

\section{Proof of EPCI in Theorem \ref{thm1}, \ref{thm2} and \ref{thm3}}

We now use Theorem \ref{thm4} to give the first proof of the entropy power concavity inequality (EPCI) for the R\'enyi entropy associated with the nonlinear diffusion
 equation $(\ref{Lup})$ on compact Riemannian manifolds with $CD(K, m)$-condition
 and compact Riemannian manifolds equipped with $(K, m)$-super Ricci flows. The advantage of this proof is that it can be easily extended to  
 prove the entropy power concavity inequality (EPCI) for the R\'enyi entropy associated with the nonlinear diffusion equation on 
 the so-called RCD$(K, N)$ metric measure spaces.  As Theorem \ref{thm1} can be regarded as the special case of Theorem \ref{thm2} for $m=n$ 
 and $\phi=0$, we need only to prove Theorem \ref{thm2}  and Theorem \ref{thm3}. 
 In Section $6$, we will give the second proof of EPCI in Theorem \ref{thm1}, \ref{thm2} and \ref{thm3}. 
 The rigidity part of Theorem \ref{thm2}  and Theorem \ref{thm3} will be also proved in Section 6.

\medskip

\begin{lemma}\label{lem1} Under the same condition as in Theorem \ref{thm4}, we have
\begin{eqnarray}
{d^2\over dt^2}E(u)&\geq &2\left(p-1+{1\over m}\right)\left( \int_M  |\nabla e'(u)|^2  u d\mu\right)^2 \left(\int_M u^p d\mu\right)^{-1}\nonumber\\
& &+2\int_M Ric_{m, n}(L)(\nabla e'(u), \nabla e'(u))u^p d\mu+\int_M {\partial g\over \partial t}(\nabla e'(u), \nabla e'(u) ud\mu. \label{E2mK1}
\end{eqnarray}
\end{lemma}
{\it Proof}. By the trace inequality we have

\begin{eqnarray*}
|\nabla^2 e'(u)|^2\geq {|\Delta e'(u)|^2\over n}.
\end{eqnarray*}
Applying the elementary inequality
\begin{eqnarray*}
(a+b)^2 \geq {a^2\over 1+\varepsilon}-{b^2\over \varepsilon}
\end{eqnarray*}
to $a=Le'(u)$, $b=\nabla\phi\cdot\nabla e'(u)$ and $\varepsilon={m-n\over n}$,  we have
\begin{eqnarray}
|\nabla^2 e'(u)|^2\geq {|L e'(u)|^2\over m}-{\nabla \phi\otimes \nabla \phi\over m-1}(\nabla e'(u), \nabla e'(u)).  \label{emK}
\end{eqnarray}
This yields
\begin{eqnarray}
{1\over 2}{d^2\over dt^2}E(u)&\geq& \int_M \left[\left(p-1+{1\over m}\right)|Le'(u)|^2+Ric_{m, n}(L)(\nabla e'(u), \nabla e'(u))\right]u^pd\mu\nonumber\\
& &\hskip2cm +{1\over 2}\int_M {\partial g\over \partial t}(\nabla e'(u), \nabla e'(u)) ud\mu.  \label{E2mK2}
\end{eqnarray}
By integration by parts and a direct calculation, we have
\begin{eqnarray}
-\int_M Le'(u) u^p d\mu&=&\int_M \nabla e'(u)\cdot \nabla u^p d\mu\nonumber\\
&=&p \int_M \nabla e'(u)\cdot u^{p-1} \nabla u d\mu\nonumber\\
&=&{p\over p-1}\int_M \nabla e'(u)\cdot \nabla u^{p-1} u d\mu\nonumber\\
&=&\int_M |\nabla e'(u) |^2 ud\mu. \label{E2mK3}
\end{eqnarray}
By the Cauchy-Schwarz inequality
\begin{eqnarray}
E'(u)^2&=&\left(\int_M |\nabla e'(u)|^2 u d\mu\right)^2=\left(\int_M (-Le'(u) )u^p d\mu\right)^2\nonumber\\
&\leq &\left(\int_M u^pd\mu\right)\left(\int_M |Le'(u)|^2 u^p d\mu \right).\label{E2mK4}
\end{eqnarray}
Combining $(\ref{E2mK2})$, $(\ref{E2mK3})$ and $(\ref{E2mK4})$, we derive $(\ref{E2mK1})$. \hfill $\square$

\medskip

\noindent{\bf Proof of EPCI in Theorem \ref{thm2}}.  In the case the metric and potential on  $M$ are time independent, under the assumption  $Ric_{m, n}(L)\geq Kg$, 
\begin{eqnarray*}
{1\over 2}{d^2\over dt^2}E(u)&\geq&\int_M \left[\left(p-1+{1\over m}\right)|Le'(u)|^2+Ric_{m, n}(L)(\nabla e'(u), \nabla e'(u))\right]u^pd\mu\\
&\geq& \left(p-1+{1\over m}\right)\int_M |Le'(u)|^2u^pd\mu+K\int_M |\nabla e'(u)|^2u^pd\mu.
\end{eqnarray*}
Combining this with $(\ref{E2mK4})$, we have

\begin{eqnarray*}
{1\over 2}{d^2\over dt^2}E(u) &\geq &\left(p-1+{1\over m}\right)\left( \int_M  |\nabla e'(u)|^2  u d\mu\right)^2\left(\int_M u^p d\mu\right)^{-1} \\
& &+K \int_M |\nabla e'(u)|^2u^p d\mu,
\end{eqnarray*}
whence
\begin{eqnarray*}
(p-1)E(u) E''(u) \geq 2\left(p-1+{1\over m}\right)E'(u)^2+2(p-1)KE(u) \left(\int_M |\nabla e'(u)|^2 u^p d\mu\right).
\end{eqnarray*}

Note that
\begin{eqnarray}
{d\over dt}N_p(u)&=&\sigma H'_pN_p(u), \nonumber\label{NHp1}\\
{d^2\over dt^2}N_p(u)&=&\left(\sigma H_p''+\sigma^2H_p'^2\right)N_p(u). \label{NHp2}
\end{eqnarray}
Now
\begin{eqnarray}
H_p'(u)={1\over p-1}{E'(u)\over E(u)}=I_p(u), \label{HIp}
\end{eqnarray}
and
\begin{eqnarray*}
H_p''(u)&=&{1\over p-1}{\partial_t E'(u)E(u)-E'(u)\partial_t E(u)\over E^2(u)}\\
&=&{1\over p-1}{E'(u)^2-E''(u)E(u)\over E^2(u)}.
\end{eqnarray*}
Hence
\begin{eqnarray*}
\sigma H_p''+\sigma^2H_p'^2&=&{\sigma\over p-1}{E'(u)^2-E''(u)E(u)\over E(u)^2}+{\sigma^2\over (p-1)^2}{E'(u)^2\over E(u)^2}\\
&=&{\sigma^2\over (p-1)^2}\left[{p-1\over \sigma}  {E'(u)^2-E''(u)E(u)\over E(u)^2}+  {E'(u)^2\over E(u)^2}\right]\\
&=&{\sigma\over (p-1)^2 E(u)^2}\left[(p-1+\sigma)  E'(u)^2- (p-1) E''(u)E(u)\right].
\end{eqnarray*}
Replacing $\sigma={2\over m}+p-1$ into the above identity, we then derive
\begin{eqnarray*}
Z:&=&(p-1+\sigma)  E'(u)^2- (p-1) E''(u)E(u)\\
&=&2\left(p-1+{1\over m}\right) E'(u)^2- (p-1) E''(u)E(u)\\
&\leq& -2(p-1)KE(u) \left(\int_M |\nabla e'(u)|^2 u^p d\mu\right).
\end{eqnarray*}
Then 
\begin{eqnarray*}
\sigma H_p''+\sigma^2H_p'^2={\sigma Z\over (p-1)^2E(u)^2}\leq -{2\sigma K \int_M |\nabla e'(u)|^2 u^p d\mu\over (p-1)E(u)}.
\end{eqnarray*}
By $(\ref{NHp2})$, we then derive  
\begin{eqnarray*}
{d^2\over dt^2}N_p(u)\leq -{2\sigma K N_p(u)\over \|u\|_p^p}\int_M |\nabla e'(u)|^2u^p d\mu.
\end{eqnarray*}
The proof of  EPCI in Theorem \ref{thm2} is completed.  \hfill $\square$

\bigskip

\noindent{\bf Proof of EPCI in Theorem \ref{thm3}}.  By Lemma \ref{lem1}, under the condition of Theorem \ref{thm3}, 
we have
\begin{eqnarray*}
{1\over 2}{d^2\over dt^2}E(u) &\geq &\left(p-1+{1\over m}\right)\left( \int_M  |\nabla e'(u)|^2  u d\mu\right)^2\left(\int_M u^p d\mu\right)^{-1} \\
& &+K_1 \int_M |\nabla e'(u)|^2u^p d\mu+K_2\int_M |\nabla e'(u)|^2ud\mu.
\end{eqnarray*}
That is to say
\begin{eqnarray*}
(p-1)E(u) E''(u) &\geq &2\left(p-1+{1\over m}\right)E'(u)^2+2(p-1)K_1E(u) \left(\int_M |\nabla e'(u)|^2 u^p d\mu\right)\\
& &\hskip2cm +2(p-1)K_2 E(u)E'(u).
\end{eqnarray*}
We then derive that
\begin{eqnarray*}
Z:&=&(p-1+\sigma)  E'(u)^2- (p-1) E''(u)E(u)\\
&=&2\left(p-1+{1\over m}\right) E'(u)^2- (p-1) E''(u)E(u)\\
&\leq& -2(p-1)K_1E(u) \left(\int_M |\nabla e'(u)|^2 u^p d\mu\right)-2(p-1)K_2E(u)E'(u),
\end{eqnarray*}
and hence
\begin{eqnarray*}
\sigma H_p''+\sigma^2H_p'^2&=&{\sigma Z\over (p-1)^2E(u)^2}\\
&\leq& -{2\sigma \over (p-1)E(u)}\left(K_1\int_M |\nabla e'(u)|^2 u^pd\mu+K_2\int_M |\nabla e'(u)|^2 ud\mu\right).
\end{eqnarray*}
By $(\ref{NHp2})$, we then derive 
\begin{eqnarray*}
{d^2\over dt^2}N_p(u)\leq -{2\sigma N_p(u) \over \|u\|_p^p}\left(K_1\int_M |\nabla e'(u)|^2 u^pd\mu+K_2\int_M |\nabla e'(u)|^2 ud\mu\right).
\end{eqnarray*}
The proof of  EPCI in Theorem \ref{thm3} is completed.  \hfill $\square$

\section{The second order derivative of Renyi entropy power}

In this section we prove an explicit formula for the second order derivative of the Renyi entropy power for the nonlinear diffusion equation $(\ref{Lup})$ 
on compact manifolds with time dependent  metrics.  When $p=1$, such an explicit formula for the second order derivative of the Shannon 
entropy power for the linear heat equation on manifolds with time dependent metrics was first proved in our previous paper \cite{LL20}.

\begin{theorem} \label{thm5}  Under the same condition and notation as in Theorem \ref{thm3}, we have 
\begin{eqnarray*}
{d^2N_p\over dt^2}
&=&-2\sigma N_p \left(  p-1+{1\over m} \right) \int_M \left[ Le'(u)-\int_M Le'(u) d\gamma \right]^2d\gamma\\
& &-2\sigma N_p\int_M Ric_{m, n}(L)(\nabla e'(u), \nabla e'(u) d\gamma-\sigma N_p\int_M {\partial g\over \partial t}(\nabla e'(u), \nabla e'(u)){u\over \|u\|_p^p}d\mu\\
& &-2\sigma N_p \int_M \left[ {m-n\over mn}\left(Le'(u)+{m\over m-n}\nabla\phi\cdot\nabla e'(u)\right)^2 + \left\|\nabla^2 e'(u)-{\Delta e'(u)\over n} g\right\|^2_{\rm HS} 
\right]d\gamma,
\end{eqnarray*}
where $d\gamma={u^pd\mu\over \|u\|_p^p}$ with  $\|u\|_p^p=\int_M u^p d\mu$. In particular, when $g$ and $\phi$ are time independent, the above result holds with 
${\partial g\over \partial t}=0$. 
\end{theorem}

In particular, when $m=n$, $\phi=0$, $L=\Delta$, $\sigma=p-1+{2\over n}$, and $g$ is time independent, we have

\begin{theorem} \label{thm5b} Under the same condition and notation as in Theorem \ref{thm1}, we have
\begin{eqnarray*}
{d^2N_p\over dt^2}
&=&-2\sigma N_p \left(  p-1+{1\over n} \right) \int_M \left[ \Delta e'(u)-\int_M \Delta e'(u) d\gamma \right]^2d\gamma\\
& &-2\sigma N_p\int_M Ric(\nabla e'(u), \nabla e'(u) d\gamma-2\sigma N_p \int_M  \left\|\nabla^2 e'(u)-{\Delta e'(u)\over n} g\right\|^2_{\rm HS} d\gamma.
\end{eqnarray*}
where $d\gamma={u^pdv\over \|u\|_p^p}$ with  $\|u\|_p^p=\int_M u^p dv$. 
\end{theorem}

We need only to prove Theorem \ref{thm5}. 
\medskip

\noindent{\bf Proof of Theorem \ref{thm5}}. Using
\begin{eqnarray*}
\|A\|_{\rm HS}^2={|{\rm tr A}|^2\over n}+\left\|A-{{\rm tr}A\over n} g\right\|^2_{\rm HS},
\end{eqnarray*}
we have

\begin{eqnarray*}
\|\nabla^2 e'(u)\|_{\rm HS}^2= {|\Delta e'(u)|^2\over n}+\left\|\nabla^2 e'(u)-{\Delta e'(u)\over n} g\right\|^2_{\rm HS}
\end{eqnarray*}
Applying the elementary equality
\begin{eqnarray*}
(a+b)^2={a^2\over 1+\varepsilon}-{b^2\over \varepsilon}+{\varepsilon\over 1+\varepsilon}\left(a+{1+\varepsilon\over \varepsilon}b\right)^2
\end{eqnarray*}
to $a=Le'(u)$, $b=\nabla\phi\cdot\nabla e'(u)$ and $\varepsilon={m-n\over n}$, we have
\begin{eqnarray*}
|\Delta e'(u)|^2= {n\over m}|L e'(u)|^2-{n\over m-n}|\nabla \phi\cdot \nabla e'(u)|^2+{m-n\over m}\left(Le'(u)+{m\over m-n}\nabla\phi\cdot\nabla e'(u)\right)^2.
\end{eqnarray*}
Therefore
\begin{eqnarray*}
\|\nabla^2 e'(u)\|_{\rm HS}^2&=&  {|L e'(u)|^2\over m}-{|\nabla \phi\cdot \nabla e'(u)|^2\over m-n}+{m-n\over mn}\left(Le'(u)+{m\over m-n}\nabla\phi\cdot\nabla e'(u)\right)^2\\
& &\hskip2cm +\left\|\nabla^2 e'(u)-{\Delta e'(u)\over n} g\right\|^2_{\rm HS}.
\end{eqnarray*}
This yields
\begin{eqnarray*}
{1\over 2}{d^2\over dt^2}E(u)&=&\int_M \left(|\nabla^2e'(u)|^2+Ric(L)(\nabla e'(u), \nabla e'(u))+(p-1)(Le'(u))^2\right)u^pd\mu\\
& &\hskip2cm +\int_M {1\over 2}{\partial g\over \partial t}(\nabla e'(u), \nabla e'(u))ud\mu\\
&=& \int_M \left[\left(p-1+{1\over m}\right)|Le'(u)|^2+Ric_{m, n}(L)(\nabla e'(u), \nabla e'(u))\right]u^pd\mu\\
& &+\int_M \left[{m-n\over mn}\left|Le'(u)+{m\over m-n}\nabla\phi\cdot\nabla e'(u)\right|^2 + \left\|\nabla^2 e'(u)-{\Delta e'(u)\over n} g\right\|^2_{\rm HS} \right]u^pd\mu\\
& &\hskip2cm + {1\over 2}\int_M {\partial g\over \partial t}(\nabla e'(u), \nabla e'(u))ud\mu.
\end{eqnarray*}
Therefore 
\begin{eqnarray*}
H_p''+\sigma H_p'^2
&=&{1\over (p-1)^2 E(u)^2}\left[(p-1+\sigma)  E'(u)^2- (p-1) E''(u)E(u)\right]\\
&=&(p-1+\sigma)\left({E'(u)\over (p-1)E(u)}\right)^2-{E''(u)\over (p-1)E(u)}\\
&=&2 \left(  p-1+{1\over m} \right) \left[\left(  \int_M Le'(u) d\gamma \right)^2-   \int_M  |Le'(u)|^2 d\gamma \right]\\
& &-2\int_M Ric_{m, n}(L)(\nabla e'(u), \nabla e'(u)) d\gamma-\int_M {\partial g\over \partial t}(\nabla e'(u), \nabla e'(u)) {u\over \|u\|_p^p}d\mu\\
& &-2\int_M \left[ {m-n\over mn}\left(Le'(u)+{m\over m-n}\nabla\phi\cdot\nabla e'(u)\right)^2 + \left\|\nabla^2 e'(u)-{\Delta e'(u)\over n} g\right\|^2_{\rm HS} \right]d\gamma.
\end{eqnarray*}
Using $(\ref{NHp1})$ and $(\ref{NHp2})$, we complete the proof of Theorem \ref{thm5}. \hfill $\square$

\section{Rigidity theorems for Renyi entropy power}

Based on Theorem \ref{thm5}, we now give the second proof of EPCI in Theorem \ref{thm1}, \ref{thm2} and \ref{thm3}. 
Moreover, we prove the rigidity part of 
 Theorem \ref{thm1}, \ref{thm2}  and \ref{thm3}.

\bigskip

\noindent{\bf Proof of Theorem \ref{thm2}}. By Theorem \ref{thm5}, under the condition $Ric_{m, n}(L)\geq Kg$, we have
\begin{eqnarray*}
{d^2N_p\over dt^2}&\leq&-{2\sigma K N_p}\int_M |\nabla e'(u)|^2 d\gamma-{2\sigma N_p}\left(  p-1+{1\over m} \right) \int_M \left[
 Le'(u)-\int_M Le'(u) d\gamma \right]^2d\gamma \\
& &\ -{2\sigma N_p}\int_M \left[ {m-n\over mn}\left(Le'(u)+{m\over m-n}\nabla\phi\cdot\nabla e'(u)\right)^2 + \left\|\nabla^2 e'(u)-{\Delta e'(u)\over n} g\right\|^2_{\rm HS} \right]d\gamma.
\end{eqnarray*}
Thus, for $p>1-{1\over m}$, we have
\begin{eqnarray*}
{d^2N_p\over dt^2}\leq-{2\sigma K N_p}\int_M |\nabla e'(u)|^2 d\gamma.
\end{eqnarray*}
Moreover,  the equality holds if and only if 
\begin{eqnarray*}
Ric_{m, n}(L)=Kg,
\end{eqnarray*}
and 
\begin{eqnarray*}
Le'(u)=\int_M Le'(u)d\gamma, \ \  (m-n)\left(Le'(u)+{m\over m-n}\nabla\phi\cdot\nabla e'(u)\right)=0, \ \ \ \nabla^2 e'(u)={\Delta e'(u)\over n} g.
\end{eqnarray*}
This yields $(\ref{rrrr1})$ in Theorem \ref{thm2}. The inequality $(\ref{Ht2})$ follows from the   Aronson-Benilan inequality $(\ref{wAB2})$  in the next section. The proof of Theorem 
\ref{thm2} is completed. 
\hfill $\square$

\bigskip

\noindent{\bf Proof of Theorem \ref{thm1}}.  Theorem \ref{thm1} follows from Theorem \ref{thm2} by taking $m=n$, $\phi=0$ and $L=\Delta$. \hfill $\square$

\bigskip

\noindent{\bf Proof of Theorem \ref{thm3}}. By Theorem \ref{thm5}, under the assumption $Ric_{m, n}(L)\geq K_1g$ and 
$\partial_t g\geq 2K_2g$,  we have
\begin{eqnarray*}
{d^2N_p\over dt^2}&\leq&-{2\sigma K_1 N_p}\int_M |\nabla e'(u)|^2 d\gamma-{2\sigma K_2 N_p\over \|u\|_p^p}\int_M |\nabla e'(u)|^2 ud\mu\\
& &-{2\sigma N_p}\left(  p-1+{1\over m} \right) \int_M \left[ Le'(u)-\int_M Le'(u)|^2 d\gamma \right]^2d\gamma \\
& &\ -{2\sigma N_p}\int_M \left[ {m-n\over mn}\left(Le'(u)+{m\over m-n}\nabla\phi\cdot\nabla e'(u)\right)^2 + \left\|\nabla^2 e'(u)-{\Delta e'(u)\over n} g\right\|^2_{\rm HS} \right]d\gamma.
\end{eqnarray*}
Hence, for $p>1-{1\over m}$, we have
\begin{eqnarray}
{d^2N_p\over dt^2}\leq-{2\sigma K_1 N_p\over \|u\|_p^p} \int_M |\nabla e'(u)|^2 u^p d\mu-{2\sigma K_2 N_p\over \|u\|_p^p}\int_M |\nabla e'(u)|^2 ud\mu.  \label{NKe}
\end{eqnarray}
Moreover,  the equality holds if and only if 
\begin{eqnarray*}
Ric_{m, n}(L)=K_1g, \ \ \ \ \partial_t g=2K_2 g, \ \ \ \partial_t \phi=nK_2,
\end{eqnarray*}
and 
\begin{eqnarray*}
Le'(u)=\int_M Le'(u)d\gamma, \ \  (m-n)\left(Le'(u)+{m\over m-n}\nabla\phi\cdot\nabla e'(u)\right)=0, \ \ \ \nabla^2 e'(u)={\Delta e'(u)\over n} g.
\end{eqnarray*}
This concludes Theorem \ref{thm3}.  \hfill $\square$

\section{NIW formula for nonlinear diffusion equation}

In his seminal paper \cite{P1}, Perelman introduced the $W$-entropy and proved its variational 
formula along the conjugate heat equation on the Ricci flow.  More precisely, let $M$ be a compact manifold equipped with the Ricci flow 
\begin{eqnarray}
\partial_t g(t)=-2Ric_{g(t)}.\label{RF}
\end{eqnarray}
Let $u(t)={e^{-f}\over (4\pi t)^{n/2}}$ be the fundamental solution to the conjugate heat equation 
\begin{eqnarray}
\partial_t u=-\Delta u+Ru.\label{BHE}
\end{eqnarray}
Let $\tau$ be such that 
$$\partial_t \tau=-1.$$
Let 
\begin{eqnarray*}
\mathcal{H}(g, u)=-\int_M u\log u dv
\end{eqnarray*}
be the Shannon entropy.  
In \cite{P1}, Perelman introduced the following $\mathcal{F}$-function and $\mathcal{W}$-entropy
\begin{eqnarray*}
\mathcal{F}(g, u)&=&\int_M (R+|\nabla \log u|^2)udv,\\
\mathcal{W}(g, u))&=&\int_M \left(\tau (R+|\nabla f|^2)+f-n\right)u dv,
\end{eqnarray*} 
and proved the following beautiful  formulae
\begin{eqnarray*} 
{d\over dt} \mathcal{H}(g(t), u(t))&=&-\mathcal{F}(g(t), u(t)), \label{FH}\\
{d\over dt} \mathcal{F}(g(t), u(t))&=&2\int_M |Ric+\nabla^2 \log u|^2 udv, \label{Ft}\\
{d \over dt}\mathcal{W}(g(\tau), u(\tau))&=&2\tau \int_M \left|Ric+\nabla^2 f-{g\over 2\tau }\right|^2 udv. \label{Wt}
\end{eqnarray*}
In particular, the $W$-entropy is always monotone increasing along the Ricci flow, and its equilibrium state (i.e., the critical point such that ${d\over dt}\mathcal{W}
  (g(\tau), u(\tau))=0$) is the so-called shrinking Ricci soliton
  \begin{eqnarray*}
Ric+\nabla^2f ={g\over 2\tau}.
\end{eqnarray*}

In \cite{N1, N2},  Ni intrduced the $W$-entropy and proved its variational formula for the heat equation $\partial_t u=\Delta u$ 
 on Riemannian manifolds with non-negative Ricci curvature. In \cite{LX}, Li and Xu extended Ni's result to Riemannian manifolds with Ricci curvature bounded 
 from below by a uniform constant. In  \cite{Li12, Li16}, the second named  author extended the $W$-entropy formula to the heat equation associated with 
the Witten Laplacian on complete Riemannian manifolds with $CD(0, m)$-condition and proved that the rigidity model of the $W$-entropy 
is the Euclidean spaces equipped with the Gaussian solitons. In \cite{LL15, LL16, LL-AJM, LL18a, LL18b, LL19a, LL19b}, the authors of this paper extended 
the $W$-entropy formula to the heat equation 
associated with the Witten Laplacian on complete Riemannian manifolds with $CD(K, m)$-condition and on compact 
$(K, m)$-super Ricci flows. In \cite{LL16, LL18b}, the authors of this paper introduced the $W$-entropy and proved the $W$-entropy formula 
along the geodesic flow on the infinite dimensional $L^2$-Wasserstein space over a complete Riemannian manifold with $CD(0, m)$-condition. 
This recovers and improves an earlier result due to Lott and Villani \cite{Lot1, LV} on the displacement convexity of the Boltzmann entropy functional on 
complete Riemannian manifold with non-negative Ricci curvature. In \cite{KL}, Kuwada and the second named author proved the monotonicity and rigidity 
theorem of the $W$-entropy on RCD$(0, N)$ metric measure spaces. 

In \cite{LL20}, the authors proved the NIW formula between the Shannon entropy power $N$, the Fisher information $I$ and the $W$-entropy for the heat 
equatioin $\partial_t u=Lu$ associated with the Witten Laplacian on complete Riemannian manifolds. More precisely, we proved the following

\begin{theorem}\label{NIW1}  Let $M$ be a complete Riemannian manifold with natural bounded geometry condition. Let $\phi\in C^4(M)$ be such that 
$\nabla^k\phi\in C_b(M)$, $k=1, 2, 3$. Let $u$ be the fundamental solution to the heat equation $\partial_t u=Lu$ associated with the Witten Laplacian
$L=\Delta-\nabla \phi\cdot\nabla$. The following NIW formula  holds
 \begin{eqnarray}
{d^2 N\over dt^2}={2N\over m}\left[{2\over m}\left(I-{m\over 2t}\right)^2+{1\over t}{d\mathcal{W}\over dt}\right],  \label{NIW-1}
\end{eqnarray}
where $N(u)=e^{-{2\over m}H(u)}$ is the Shannon entropy power, $I(u)=\int_M {|\nabla u|^2\over u}d\mu$ is the Fisher information, and $\mathcal{W}$ is the 
$W$-entropy defined by 
\begin{eqnarray}
\mathcal{W}(u)=\int_M \left(t |\nabla f|^2+f-m\right)u d\mu.\label{Wm-1}
\end{eqnarray}
As a consequence of $(\ref{NIW-1})$, under the $CD(0, m)$-condition, we have
\begin{eqnarray}
{d\mathcal{W}(u)\over dt}\leq -{2t\over m}\left|I(u)-{m\over 2t}\right|^2. \label{NIW-2}
\end{eqnarray}
In particular, the NIW formula $(\ref{NIW-1})$ and the $W$-entropy inequality $(\ref{NIW-2})$ hold for 
the heat equation $\partial_t u=\Delta u$ with $m=n$, $\mu=v$ and $L=\Delta$  on $n$-dimensional complete  Riemannian manifolds 
with non-negative Ricci curvature and  with 
natural bounded geometry condition. 
\end{theorem}

Moreover,  we proved the NFW formula in  \cite{LL20}  for the Shannon entropy power $\mathcal{N}$, the Perelman $\mathcal{F}$-functional and the Perelman $\mathcal{W}$-entropy 
 along the conjugate heat equation on the Ricci flow. It provide us a new understanding for Perelman's mysterious $\mathcal{W}$-
entropy formula for Ricci flow and can be used to characterize the shrinking Ricci solitons. 

\begin{theorem}\label{NFW3}  Let $(M, g(t), t\in [0, T])$ be a compact Ricci flow $\partial_t g=-2Ric_{g(t)}$. Let $u(t)$ be the fundamental solution to the conjugate heat equation 
\begin{eqnarray*}
\partial_t u=-\Delta u+Ru.\label{BHE1}
\end{eqnarray*}
Let
\begin{eqnarray*}
H(u(t))=-\int_M u\log udv, \ \ \ \ \mathcal{N}(u(t))=e^{{2\over n}H(u(t))}.
\end{eqnarray*}
Then 
\begin{eqnarray*}
{d^2 \mathcal{N} \over dt^2}={2\mathcal{N} \over n}\left[{2\over n}\left(\mathcal{F}-{n\over 2\tau}\right)^2+{1\over \tau}{d\mathcal{W}\over dt}\right],  \label{NIWRF}
\end{eqnarray*}
where 
\begin{eqnarray*}
\mathcal{F}=\int_M (R+|\nabla \log u|^2)udv,
\end{eqnarray*}
and
\begin{eqnarray*}
{d\mathcal{W}\over dt}=2 \tau \int_M \left|Ric+\nabla^2 f-{g\over 2\tau }\right|^2udv.
\end{eqnarray*}
In particular,  the Shannon entropy power is convex on $(0, T]$, i.e., 
\begin{eqnarray*}
{d^2\over dt^2} \mathcal{N}(u(\tau))\geq 0.
\end{eqnarray*}
Moreover,  ${d^2\over dt^2}\mathcal{N}(u(\tau))=0$ holds at some $\tau=\tau_0\in (0, T]$  if and only if  $(M, g(\tau), u(\tau))$  is a shrinking Ricci soliton 
\begin{eqnarray*}
Ric+\nabla^2 \log u={g\over 2\tau}.
\end{eqnarray*} 
\end{theorem}

In \cite{LNVV}, Lu-Ni-Vazquez-Villani introduced the $W$-entropy for the porous medium equation and the fast diffusion equation on  compact Riemannian manifolds and proved their 
monotonicity on compact Riemannian manifolds with non-negative Ricci curvature. Inspired by our previous work on the NIW formula and NIW inequality for the heat equation on 
Riemannian manifolds,  it is natural to ask the question whether we can extend the NIW formula and the NIW  inequality to the porous medium  equation and the fast 
diffusion equation on Riemannian manifolds. The purpose of this section is to study this problem. 

\subsection{NIW formula for NLDE for Laplacian}

In this section we prove the NIW formula for nonlinear diffusion equation for usual Laplacian on compact Riemannian manifolds. To state the main result of this  subsection, let us briefly 
recall  Lu-Ni-Vazquez-Villani's work \cite{LNVV} on the entropy formulae for the nonlinear diffusion equation $\partial_t u=\Delta u^p$ on compact Riemannian manifolds.

Let $M$ be an $n$-dimensional compact Riemannian manifold. Let $u$ be a  positive and smooth solution to the nonlinear diffusion equation 

\begin{eqnarray}
\partial_t u=\Delta u^p, \label{NLDE}
\end{eqnarray}
 where $p>1-{1\over n}$.  Let $v={p\over p-1} u^{p-1}$. Then 
\begin{eqnarray*}
\partial_t v=(p-1)v\Delta v+|\nabla v|^2.
\end{eqnarray*}
For any $\alpha\geq 1$, let 
\begin{eqnarray*}
F_\alpha=(p-1)\Delta v+(\alpha-1){|\nabla v|^2\over v},
\end{eqnarray*}
in particular, $F_1=(p-1)\Delta v$. Let $\kappa={n\over  n(p-1)+2}$, $a=(p-1)\kappa$, and $b=n(p-1)$. Define 
\begin{eqnarray}
\mathcal{N}_u(t)=-t^a \int_M v udv={p\over 1-p}t^{a}\int_M u^pdv. \label{Nu-1}
\end{eqnarray}
By Lemma 5.1 in \cite{LNVV}, it holds
\begin{eqnarray}
{d\over dt}\mathcal{N}_u(t)=-t^a\int_M \left(F_1+{a\over t}\right)vudv. \label{Nu-2}
\end{eqnarray}
Note that when $Ric\geq 0$, the Aronson-Benilan inequality was proved in \cite{LNVV}: When $p>1$, it holds
\begin{eqnarray}
F_1+{a\over t}\geq 0, \label{AB-1}
\end{eqnarray}
and the inverse inequality in $(\ref{AB-1})$ holds when $p\in \left(1-{2\over n}, 1\right)$. Equivalently, for all $p>1-{2\over n}$ and $p\neq 1$, it holds
\begin{eqnarray}
\Delta v+{\kappa\over t}\geq 0. \label{AB-12}
\end{eqnarray}

Following \cite{LNVV}, we introduce the Perelman $W$-entropy associated with $(\ref{NLDE})$  as follows
\begin{eqnarray}
\mathcal{W}_u(t):=t{d\over dt}\mathcal{N}_u(t)+\mathcal{N}_u(t). \label{Wp-1}
\end{eqnarray}
To emphasis the dependence on $p$, we use the notation $\mathcal{W}_p(u)$ instead of $\mathcal{W}_u$. By \cite{LNVV}, we have
\begin{eqnarray}
\mathcal{W}_p(u)=t^{a+1}\int_M \left(p{|\nabla v|^2\over v}-{a+1\over t}\right)vudv. \label{Wp-2}
\end{eqnarray}
The following $W$-entropy formula is due to Lu-Ni-Vazquez-Villani  (Theorem 5.5 in \cite{LNVV}).
\begin{eqnarray}
{d\over dt}\mathcal{W}_p(u)&=&-2(p-1)t^{a+1}\int_M \left(\left|v_{ij}+{g_{ij}\over (b+2)t}\right|^2+R_{ij}v_iv_j\right)vudv\nonumber\\
& &\hskip2cm -2t^{a+1}
\int_M \left(F_1+{a\over t}\right)^2vudv.   \label{Wp-3}
\end{eqnarray}
 
Recall that $e(u)={u^p\over p-1}$, $e'(u)={p\over p-1}u^{p-1}$. Thus $v=e'(u)$, and we can reformulate $(\ref{Wp-3})$ as follows
\begin{eqnarray}
{1\over t^{a+1}}{d\over dt}\mathcal{W}_p(u)&=&-2p\int_M \left(\left|\nabla^2 e'(u)+{g\over (n(p-1)+2)t}\right|^2+Ric(\nabla e'(u), \nabla e'(u))\right)u^pdv \nonumber\\
& &\hskip2cm -{2p\over p-1}\int_M \left((p-1)\Delta e'(u)+{a\over t}\right)^2u^pdv.   \label{Wp-4}
\end{eqnarray}

On the other hand, by Theorem \ref{thm5}, the second order derivative of the Renyi entropy power associated with the 
nonlinear diffusion equation $(\ref{NLDE})$  is given by

\begin{eqnarray}
{\|u\|_p^p\over 2\sigma N_p(u)}{d^2 N_p(u)\over dt^2}&=&
-\left(p-1+{1\over n}\right)\int_M \left|\Delta e'(u)-\int_M \Delta e'(u) d\gamma\right|^2u^p dv\nonumber\\
& &\hskip2cm -\int_M Ric(\nabla e'(u), \nabla e'(u)) u^pdv\nonumber\\
& &\hskip2cm -\int_M \left\|\nabla^2 e'(u)-{\Delta e'(u)\over n}g\right\|_{\rm HS}^2 u^pdv.   \label{Wp-5}
\end{eqnarray}
where $d\gamma={u^p dv\over \int_M u^p dv}={vudv\over \int_M vudv}$. 

\medskip

Now we state the main result of this subsection.

\begin{theorem}\label{NIW2} Let $M$ be a compact Riemannian manifold, $u$ be a positive smooth solution to the nonlinear diffusion equation 
$(\ref{NLDE})$ with $p>0$.  Then, under the above notations,  we have the following NIW formula 
 \begin{eqnarray*}
 {d^2 N_p(u)\over dt^2}=
2\sigma N_p(u) \left(  {1+n(p-1) \over n} \left|I_p(u)-{\kappa\over t} \right|^2+{1\over 2p \|u\|_p^p}{1\over t^{a+1}}{d\over dt}\mathcal{W}_p(u)\right).
 \end{eqnarray*}Under the condition $Ric\geq 0$, we have
 \begin{eqnarray*}
 {1\over 2p  \|u\|_p^p}{1\over  t^{a+1}}
{d\over dt}\mathcal{W}_p(u)\leq  - {1+n(p-1) \over n}\left|I_p(u)-{\kappa\over t} \right|^2. 
 \end{eqnarray*}
 In particular, under the condition $Ric\geq 0$, the equality in the above inequality holds if and only if ${d^2 N_p\over dt^2}=0$, which is equivalent to say  that $(M, g)$ is Ricci flat, i.e., $Ric=0$, and $\nabla^2 v=-{I_p(u)\over n} g$, 
where $I_p(u):={d\over dt}H_{p}(u(t))$ is the $p$-th Fisher information and satisfies the differential equation 
\begin{eqnarray*}
I'_p(u)+\sigma I_p^2(u)=0.
\end{eqnarray*}
\end{theorem}
{\it Proof}. To simplify the notation, let $I_1= {\|u\|_p^p\over 2\sigma N_p(u)}{d^2 N_p(u)\over dt^2}$, $I_2={1\over 2p t^{a+1}}{d\over dt}\mathcal{W}_p(u)$, and 
$I_3=I_1-I_2$. By $(\ref{Wp-4})$ and $(\ref{Wp-5})$, we have 
\begin{eqnarray}
I_3&=&\int_M \left\|\nabla^2 e'(u)+{g\over (n(p-1)+2)t}\right\|_{\rm HS}^2u^pdv-
\int_M \left\|\nabla^2 e'(u)-{\Delta e'(u)\over n}g\right\|_{\rm HS}^2u^pdv\nonumber\\
& &+{1\over p-1}\int_M \left((p-1)\Delta e'(u)+{a\over t}\right)^2u^pdv-\left(p-1+{1\over n}\right)\int_M \left|\Delta e'(u)-\int_M \Delta e'(u) d\gamma\right|^2u^pdv\nonumber\\
&:=&I_{31}+I_{32}.
\end{eqnarray}
Now 
\begin{eqnarray*}
I_{31}&=&\int_M \left[{2\over (2+n(p-1))t}\Delta e'(u)+{1\over n}|\Delta e'(u)|^2+{n\over (n(p-1)+2)^2 t^2}\right] u^pdv,\nonumber\\
I_{32}&=&\int_M \left[-{1\over n}|\Delta e'(u)|^2 +{2a\over t}\Delta e'(u) +{a^2\over (p-1)t^2}\right]u^pdv\nonumber\\
& &\hskip1cm +\left(p-1+{1\over n}\right)
\left(\int_M \Delta e'(u) u^p dv\right)^2\left(\int_M u^p dv\right)^{-1},
\end{eqnarray*}
which yields
\begin{eqnarray*}
I_3&=&{1\over t^2} \left({n\over (2+n(p-1))^2} +{a^2\over p-1}\right)\int_M u^p dv+{2\over t}\left(a+{1\over 2+n(p-1)}\right)\int_M \Delta e'(u) u^p dv\\
& &\hskip2cm +\left(p-1+{1\over n}\right)\left(\int_M \Delta e'(u) u^p dv\right)^2\left(\int_M u^p dv\right)^{-1}\\
&=& {1\over t^2} {n(1+n(p-1))\over (2+n(p-1))^2}\int_M u^p dv+{2\over t}{1+n(p-1)\over 2+n(p-1)}\int_M \Delta e'(u) u^p dv\\
& &\hskip2cm +{1+n(p-1)\over n}\left(\int_M \Delta e'(u) u^p dv\right)^2\left(\int_M u^p dv\right)^{-1}\\
&=& {1+n(p-1)\over n}\left[\int_M \Delta e'(u) u^p dv+{n\over (2+n(p-1))t}\int_M u^pdv  \right]^2\left[\int_M u^p dv\right]^{-1}\nonumber\\
&=& {1+n(p-1)\over n}\left[\int_M \Delta e'(u) d\gamma+{n\over (2+n(p-1))t} \right]^2\left[\int_M u^p dv\right]\nonumber.
\end{eqnarray*}
Note that $\kappa={n\over 2+n(p-1)}$ and integration by parts yields
\begin{eqnarray}
\int_M \Delta e'(u) d\gamma&=&{\int_M\Delta e'(u) u^pdv\over \int_M u^p dv} =-{\int_M \nabla e'(u)\cdot \nabla u^p dv\over \int_M u^pdv}\nonumber\\
&=& -{p^2\int_M u^{2p-3}|\nabla u|^2 dv\over \int_M u^pdv}=-p^2\int_M u^{p-3} |\nabla u|^2 d\gamma\nonumber\\
&=&-I_p(u). \label{LvIp}
\end{eqnarray}
Hence
\begin{eqnarray*}
I_3= {1+n(p-1)\over n}\|u\|_p^p\left|I_p(u)-{\kappa\over t} \right|^2.
\end{eqnarray*}
That is to say
 \begin{eqnarray*}
 {\|u\|_p^p\over 2\sigma N_p(u)}{d^2 N_p(u)\over dt^2}=
 {1+n(p-1)\over n}\|u\|_p^p\left|I_p(u)-{\kappa\over t} \right|^2+{1\over 2p t^{a+1}}{d\over dt}\mathcal{W}_p(u).
 \end{eqnarray*}
 This completes the proof of Theorem \ref{NIW2}. \hfill $\square$

\begin{remark} In particular, taking limit $p\rightarrow 1$  in Theorem \ref{NIW2}, we recover the NIW formula $(\ref{NIW-1})$ with $m=n$ in Theorem \ref{NIW1} for the Shannon entropy power $N(u)$, the Fisher information $I(u)$ and the $W$-entropy associated with 
 the heat equation $\partial_t u=\Delta u$ on Riemannian manifold, i.e., 
\begin{eqnarray*}
 {d^2 N(u)\over dt^2}=
{2\over n} N(u) \left[  {2\over n}\left|I(u)-{n\over 2t} \right|^2+{1\over  t}{d\over dt}\mathcal{W}(u)\right].
 \end{eqnarray*}
 Under the condition $Ric\geq 0$, we have
  \begin{eqnarray*}
 {d\mathcal{W}(u)\over dt}\leq  - {2t \over n}\left|I(u)-{n\over 2t} \right|^2. 
 \end{eqnarray*}
In particular, under the condition $Ric\geq 0$, the equality in the above inequality holds if and only if ${d^2 N\over dt^2}=0$, which is equivalent to say  that $(M, g)$ is Ricci flat, i.e., $Ric=0$, and $v=\log u$ satisfies $\nabla^2 v=-{I(u)\over n} g$, 
where $I(u)$ is a solution to the differential equation 
\begin{eqnarray*}
I'_p(u)+{2\over n} I_p^2(u)=0.
\end{eqnarray*}

\end{remark}

 \subsection{Aronson-Benilan estimates for NLDE  for Witten Laplacian}

 Let $u$ be a positive and smooth solution to the nonlinear diffusion equation $(\ref{Lup})$ associated with the Witten Laplacian on $(M, g, \mu)$, i.e., 
\begin{eqnarray}
\partial_t u=L u^p, \label{WNLDE1}
\end{eqnarray}
 where $p>1-{1\over m}$.  
 Let $v={p\over p-1} u^{p-1}$. Then 
\begin{eqnarray*}
\partial_t v=(p-1)vL v+|\nabla v|^2.
\end{eqnarray*}
For any $\alpha\geq 1$, let 
\begin{eqnarray*}
F_\alpha=(p-1)Lv+(\alpha-1){\partial_t v\over v}=\alpha (p-1)Lv+(\alpha-1){|\nabla v|^2\over v}. \label{Fa}
\end{eqnarray*}
In particular, $F_1=(p-1)L v$. 

\begin{proposition} Let $\mathcal{L}=\partial_t-(p-1)vL$. 
Then
\begin{eqnarray*}
\mathcal{L} F_\alpha=2(p-1)\left( |\nabla^2 v|^2+Ric(L)(\nabla v, \nabla v)\right)+2p\langle \nabla F_\alpha, \nabla v\rangle +(\alpha-1) \left({\partial_t v\over v}\right)^2+F_1^2. \label{LFa}
\end{eqnarray*}
In particular, for $\alpha=1$, it holds
\begin{eqnarray}
\mathcal{L} F_1=2(p-1)\left[ |\nabla^2 v|^2+Ric(L)(\nabla v, \nabla v)\right]+2p\langle \nabla F_1, \nabla v\rangle +F_1^2.   \label{LF1}
\end{eqnarray}
\end{proposition}
{\it Proof}.  The proof is similar to the one of Proposition 3.2 in \cite{LNVV} and uses the generalized Bochner formula.  To save the length of the paper, we omit the detail. \hfill $\square$

\medskip

We  have the following Aronson-Benilan estimate for the nonlinear diffusion equation $(\ref{Lup})$ associated with the Witten Laplacian for $p>1-{2\over m}$ and $p\neq 1$, which extends the corresponding result due to 
Lu-Ni-Vazquez-Villani \cite{LNVV}. When $p=1$, it  corresponds 
the Li-Yau differential Harnack inequality for the heat equation $\partial_t u=Lu$ proved in \cite{Li05, Li12}.

\begin{theorem}\label{thmwAB} Let $u$ be a positive and smooth solution to  the nonlinear diffusion equation $(\ref{Lup})$ for $p>1-{2\over m}$ and $p\neq 1$. Suppose that $Ric_{m, n}(L)\geq 0$.  
Then the Aronson-Benilan Harnack inequality holds: when $p>1$,  
\begin{eqnarray}
F_1+{(p-1)\kappa \over t}\geq 0,  \label{wAB1a}
\end{eqnarray}
and when $p\in (1-{2\over m}, 1)$, 
\begin{eqnarray}
F_1+{(p-1)\kappa \over t}\leq 0,  \label{wAB1b}
\end{eqnarray}
Equivalently, for all $p>1-{2\over m}$, it holds 
\begin{eqnarray}
Lv+{\kappa \over t}\geq 0, \label{wAB2}
\end{eqnarray}
where 
\begin{eqnarray*}
\kappa={m\over  m(p-1)+2}.
\end{eqnarray*}

\end{theorem}
{\it Proof}. The proof is similar to the one of  the Aronson-Benilan estimate $(\ref{AB-1})$ and $(\ref{AB-12})$  for the PME or FDE associated with the Laplacian on compact Riemannian manifolds with non-negative 
Ricci curvature \cite{LNVV}. Indeed,  for any $m\geq n$, we have (cf. \cite{BE, Li05}) 
\begin{eqnarray*}
 |\nabla^2 v|^2+Ric(L)(\nabla v, \nabla v)\geq {|Lv|^2\over m}+Ric_{m, n}(L)(\nabla v, \nabla v). \label{CD1}
\end{eqnarray*}
Thus, as $Ric_{m, n}(L)\geq 0$, we have 
\begin{eqnarray*}
 |\nabla^2 v|^2+Ric(L)(\nabla v, \nabla v)\geq {|Lv|^2\over m}. \label{CD2}
\end{eqnarray*}
By  $(\ref{LF1})$,  for $p>1$,  we have
\begin{eqnarray*}
\mathcal{L} F_1\geq {2(p-1)\over m} |Lv|^2+2p\langle \nabla F_1, \nabla v\rangle +F_1^2,  \label{LF2}
\end{eqnarray*}
whence
\begin{eqnarray}
\mathcal{L} F_1\geq \left({2\over m(p-1)}+1\right) F_1^2+2p\langle \nabla F_1, \nabla v\rangle.   \label{LF3}
\end{eqnarray}
Let $F=tF_1$. By calculation and $(\ref{LF3})$, we have
\begin{eqnarray*}
\mathcal{L} F&=& t\mathcal{L} F_1+F_1\nonumber\\
&\geq& \left({2\over m(p-1)}+1\right) {F^2\over t}+2p\langle \nabla F, \nabla v\rangle+{F\over t}.  \label{LF4}
\end{eqnarray*}
from which and the maximum principle, we can prove 
\begin{eqnarray*}
F\geq -(p-1)\kappa.
\end{eqnarray*}
 Equivalently, $(\ref{wAB1a})$ or $(\ref{wAB2})$ holds.  
 
 On the other hand, by  $(\ref{LF1})$, for $p\in (1-{2\over m}, 1)$,  we have
 \begin{eqnarray}
\mathcal{L} F_1-2p\langle \nabla F_1, \nabla v\rangle\leq F_1^2-{2(1-p)\over m}|Lv|^2=-\left({2\over m(1-p)}-1\right) F_1^2.   \label{LF5}
\end{eqnarray}
 Let $F=tF_1$. By calculation and $(\ref{LF5})$, we have
 \begin{eqnarray*}
\mathcal{L} F&=& t\mathcal{L} F_1+F_1\nonumber\\
&\leq& -\left({2\over m(1-p)}-1\right) {F^2\over t}+2p\langle \nabla F, \nabla v\rangle+{F\over t}.  \label{LF6}
\end{eqnarray*}
 from which and the maximum principle, we can prove 
\begin{eqnarray*}
F\leq -(p-1)\kappa.
\end{eqnarray*}
 Equivalently, $(\ref{wAB1b})$ or $(\ref{wAB2})$ holds.  The proof of Theorem \ref{thmwAB} is completed. 
  \hfill $\square$ 

\medskip

 \subsection{$W$-entropy formula for NLDE for Witten Laplacian}
 
 In this subsection we extend the $W$-entropy formula due to  Lu-Ni-Vazquez-Villani \cite{LNVV}  to nonlinear diffusion  equation associated with 
 the Witten Laplacian on compact Riemannian manifolds. See related work by Huang and Li \cite{HL}.

 Let $(M, g)$ be a compact Riemannian manifold, $m\geq n$, $\kappa={m\over  m(p-1)+2}$,  $a=(p-1)\kappa$, and $b=m(p-1)$.  Let $u$ be a positive and smooth solution to the nonlinear diffusion equation 
$ (\ref{WNLDE1})$  associated with the Witten Laplacian on $(M, g, \mu)$. Inspired by Lu-Ni-Vazquez-Villani \cite{LNVV}, define 
\begin{eqnarray*}
\mathcal{N}_u(t)=-t^a \int_M v ud\mu. \label{WNu-1}
\end{eqnarray*}
Similarly to Lemma 5.1, Lemma 5.2   and $(5.4)$  in \cite{LNVV}, we have

\begin{lemma}\label{lemm1} Let $p>0$, and let $v$  and $F_1$ be as in Theorem\ref{thmwAB}. Then
\begin{eqnarray*}
{d\over dt}\int_M v ud\mu &=& \int_M F_1 vu d\mu=-p\int_M |\nabla v|^2 ud\mu, \\
{d\over dt}\int_M F_1 vu d\mu &=& 2\int_M \left((p-1)(|\nabla^2 v|^2+Ric(L)(\nabla v, \nabla v)+ F_1^2\right)vud\mu.
\end{eqnarray*}
Moreover
\begin{eqnarray*}
{d\over dt}\mathcal{N}_u(t)=-t^a\int_M \left(F_1+{a\over t}\right)vud\mu. 
\end{eqnarray*}
\end{lemma}
{\it Proof}. The proof is similar to the one in \cite{LNVV}. \hfill $\square$

\medskip

Following \cite{LNVV}, we introduce the Perelman $W$-entropy associated with $(\ref{NLDE})$  as follows
\begin{eqnarray}
\mathcal{W}_u(t):=t{d\over dt}\mathcal{N}_u(t)+\mathcal{N}_u(t). \label{WWp-1}
\end{eqnarray}
To emphasis the dependence on $p$, we use the notation $\mathcal{W}_p(u)$ instead of $\mathcal{W}_u$. By \cite{LNVV}, we have
\begin{eqnarray*}
\mathcal{W}_p(u)=t^{a+1}\int_M \left(p{|\nabla v|^2\over v}-{a+1\over t}\right)vud\mu. \label{WWp-2}
\end{eqnarray*}

\begin{theorem}\label{thmWp} Let $m\geq n$, and let $u$ be a positive smooth solution to the 
nonlinear diffusion equation $(\ref{WNLDE1})$ with $p>0$. Let $v={p\over p-1} u^{p-1}$, $F_1=(p-1)L v$, 
$\kappa={m\over  m(p-1)+2}$, $a=(p-1)\kappa$ and $b=m(p-1)$. Define the $W$-entropy by $(\ref{WWp-1})$. Then 
\begin{eqnarray}
{d\over dt}\mathcal{N}_u(t)=-t^a\int_M \left(F_1+{a\over t}\right)vud\mu. \label{WNu-2}
\end{eqnarray}
and
 \begin{eqnarray}
{d\over dt}\mathcal{W}_p(u)&=&-2pt^{a+1}\int_M \left[\left\|\nabla^2v+{g\over (b+2)t}\right\|_{\rm HS}^2 +
{1\over m-n}\left(\nabla \phi\cdot\nabla v-{m-n\over (b+2)t}\right)^2\right]u^pd\mu\nonumber\\
& &-2pt^{a+1}\int_M Ric_{m, n}(L)(\nabla v, \nabla v)u^pd\mu-2p(p-1)t^{a+1}
\int_M \left(Lv+{\kappa\over t}\right)^2u^pd\mu. \nonumber\\
& & \label{WWp-4}
\end{eqnarray}
\end{theorem}
{\it Proof}. We can combine the argument used in the proof of the $W$-entropy formula for $\partial_t u=\Delta u^p$ in Theorem 5.5 in  Lu-Ni-Vazquez-Villani \cite{LNVV}  
and the argument used in the proof of the $W$-entropy formula   for $\partial_t u=Lu$ in Theorem 2.3 in \cite{Li12} to give a proof of Theorem \ref{thmWp}. To save the length of the paper, 
we omit the detail here.   \hfill $\square$

Note that, there are two terms in the integrant in the right hand side of the first line of the $W$-entropy formula $(\ref{WWp-4})$. To better understand the geometric meaning of the sum of these 
two terms, we would like to give an alternative proof of the $W$-entropy formula $(\ref{WWp-4})$ for the case $m\in \mathbb{N}$ with $m\geq n$. In this case, let $(N, g_N)$ be a $(m-n)$ dimensional compact
 Riemannian manifold with volume $v_N(N)=1$, and let  $\widetilde{M}=M\times N$ be equipped with the warped product metric
 \begin{eqnarray*}
 \widetilde{g}=g\oplus e^{-{2\phi\over m-n}}g_N. 
 \end{eqnarray*}
 Then the volume measure on $\widetilde{M}$ is 
 \begin{eqnarray*}
 dvol_{\widetilde{M}}=d\mu\otimes dv_N.
 \end{eqnarray*}
 
 Recall the following 
 \begin{proposition} (\cite{LL15})\label{LD}
 Let $\Delta_{\widetilde{M}}$ be the Laplace-Beltrami operator on $(\widetilde{M}, \widetilde{g})$. Then 
 \begin{eqnarray*}
 \Delta_{\widetilde{M}} =L+e^{-{2\phi\over m-n}}\Delta_N,
 \end{eqnarray*}
 where $\Delta_N$ is the Laplace-Beltrami operator on $(N, g_N)$. In particular, for   $f\in C^2(M)$, 
  \begin{eqnarray*}
 \Delta_{\widetilde{M}} f =Lf.
 \end{eqnarray*}
 Moreover, for  $f\in C^2(M)$, we have
 \begin{eqnarray*}
 \widetilde{\nabla}_if=\nabla_i f, \ \ \ \widetilde{\nabla}_\alpha f=0,
 \end{eqnarray*}
 and
\begin{eqnarray*}
\widetilde{\nabla}^2_{ij} f={\nabla}^2_{ij} f, \ \ \ \widetilde{\nabla}^2_{i\alpha} f=0, \ \ \ \widetilde{\nabla}^2_{\alpha\beta} f=-{\nabla \phi\cdot\nabla f\over m-n}g_{\alpha\beta}.
\end{eqnarray*}
where $\widetilde{\nabla}$ and $\widetilde{\nabla}^2$ denote the gradient and the Hessian on $(\widetilde{M}, \widetilde{g})$, 
$i, j=1, \ldots, n$ denotes the direction of local coordinates along $M$, and $\alpha, \beta=n+1, \ldots, m$ denotes the 
direction of local coordinates along $N$. 
 \end{proposition}

 In view of Proposition \ref{LD}, the nonlinear diffusion equation $(\ref{WNLDE1})$  associated with the Witten Laplacian  is equivalent to  the 
 nonlinear diffusion equation associated with the Laplace-Beltrami operator on $\Delta_{\widetilde{M}}$ on 
 $(\widetilde{M}, \widetilde{g})$, i.e.,  
\begin{eqnarray}
\partial_t u= \Delta_{\widetilde{M}}u^p, \label{WNLDE2}
\end{eqnarray}

By the $W$-entropy formula $(\ref{Wp-3})$  due to  Lu-Ni-Vazquez-Villani \cite{LNVV} for the nonlinear diffusion equation $(\ref{WNLDE2})$ associated with the Laplace-Beltrami operator 
on $(\widetilde{M}, \widetilde{g})$, we have
\begin{eqnarray}
{d\over dt}\mathcal{W}_p(u)&=&-2(p-1)t^{a+1}\int_M \left(\left\|\widetilde{\nabla}^2v+{\widetilde{g}\over (b+2)t}\right\|_{\rm HS}^2+\widetilde{Ric}
(\widetilde{\nabla} v, \widetilde{\nabla} v)\right)vud\mu\nonumber\\
& &\hskip2cm -2t^{a+1}
\int_M \left(F_1+{a\over t}\right)^2vud\mu,  \label{WWp-3}
\end{eqnarray}
 where $\widetilde{Ric}$ denotes the Ricci curvature on $(\widetilde{M}, \widetilde{g})$. By \cite{Bess, Lot1, Li05, LL15}, we have
  \begin{eqnarray*}\widetilde{Ric}
(\widetilde{\nabla} v, \widetilde{\nabla} v)={Ric}_{m, n}(L)(\nabla v, \nabla v). \label{WRic}
 \end{eqnarray*}
 On the other hand, by Proposition \ref{LD}, we have 
 \begin{eqnarray}
\left\|\widetilde{\nabla}^2v+{\widetilde{g}\over (b+2)t}\right\|_{\rm HS}^2
=\left\|\nabla^2v+{g\over (b+2)t}\right\|_{\rm HS}^2+{1\over m-n}\left(\nabla \phi\cdot\nabla v-{m-n\over (b+2)t}\right)^2. \label{WHess1}
\end{eqnarray}
Therefore, we reprove the $W$-entropy formula for the  nonlinear diffusion equation $(\ref{WNLDE1})$ 
associated with the Witten Laplacian on $(M, g, \mu)$ for the case $m\in \mathbb{N}$ and $m\geq n$. 

\begin{remark} \label{Wpm}
The above proof is an analogue of the one  in our previous paper \cite{LL15} for the $W$-entropy formula for
 the linear heat equation $\partial_t u=Lu$ associated with the Witten Laplacian on $(M, g, \mu)$ for $m\in \mathbb{N}$ 
 with $m\geq n$.  The advantage of this proof is that it gives a geometric interpretation of the integrant in the first line of the right had side of the $W$-entropy formula. Namely, by $(\ref{WHess1})$, the 
 quantity $\left\|\nabla^2v+{g\over (b+2)t}\right\|_{\rm HS}^2 +
{1\over m-n}\left(\nabla \phi\cdot\nabla v-{m-n\over (b+2)t}\right)^2$ equal to $\left\|\widetilde{\nabla}^2v+{
\widetilde{g}\over (b+2)t}\right\|_{\rm HS}^2$ which appear in the first integrant in the right hand side of 
Lu-Ni-Vazquez-Villani's $W$-entropy formula  $(\ref{WWp-3})$ for the  
nonlinear diffusion equation $(\ref{WNLDE2})$ on the warped product space  $(\widetilde{M}, \widetilde{g})$. 
\end{remark}

\begin{corollary}\label{corWp2}  Let $(M, g)$ be a compact Riemannian manifold with $CD(0, m)$-condition (i.e., $Ric_{m, n}(L)\geq 0$). Let $u$ be a 
positive and smooth solution to  the 
nonlinear diffusion equation $(\ref{WNLDE1})$. Then \\
$(i)$ For $p>1-{2\over m}$, ${d\over dt}\mathcal{N}_u(t)\leq 0$, and any ancient positive smooth solution to $(\ref{WNLDE1})$ must be a constant. \\
$(ii)$ For $p\geq 1-{1\over m}$,  ${d\over dt}\mathcal{W}_p(t)\leq 0$. In particular, $\mathcal{N}_u(t)$ is a monotone non-decreasing concave function in ${1\over t}$. 
\end{corollary} 
{\it Proof}. The proof is similar to the one in
\cite{LNVV}  for the nonlinear diffusion equation $\partial_t u=\Delta u^p$ on compact Riemannian manifold with  negative Ricci curvature.  \hfill $\square$

 \subsection{NIW formula for NLDE for Witten Laplacian}
 
 In this subsection we prove the NIW formula for nonlinear diffusion  equation associated with 
 the Witten Laplacian on compact Riemannian manifolds.

\begin{theorem}\label{NIW3} Let $m\geq n$,  and $u$ be a positive and  smooth solution to the 
nonlinear diffusion equation $(\ref{WNLDE1})$. Let $v={p\over p-1} u^{p-1}$, $F_1=(p-1)L v$, 
$\kappa={m\over  m(p-1)+2}$, $a=(p-1)\kappa$ and $b=m(p-1)$. Define the $W$-entropy by $(\ref{WWp-1})$.
Then the following NIW formula holds
 \begin{eqnarray*}
 {d^2 N_p(u)\over dt^2}=
2\sigma N_p(u) \left(  {1+m(p-1) \over m}\left|I_p(u)-{\kappa\over t} \right|^2+{1\over 2p  \|u\|_p^p}{1\over  t^{a+1}}
{d\over dt}\mathcal{W}_p(u)\right). \label{NIW4}
 \end{eqnarray*}
 Under the condition $CD(0, m)$, we have
 \begin{eqnarray*}
 {1\over 2p  \|u\|_p^p}{1\over  t^{a+1}}
{d\over dt}\mathcal{W}_p(u)\leq  - {1+m(p-1) \over m}\left|I_p(u)-{\kappa\over t} \right|^2. \label{NIW4b}
 \end{eqnarray*}
 In particular, under the condition $Ric_{m, n}(L)\geq 0$, the equality in the above inequality holds if and only if ${d^2 N_p\over dt^2}=0$, which is equivalent to say  that $(M, g)$ is quasi-Ricci flat, i.e., $Ric_{m, n}(L)=0$, and $\nabla^2 v=-{I_p(u)\over m} g$, 
where $I_p(u):={d\over dt}H_{p}(u(t))$ is the $p$-th Fisher information and satisfies the differential equation 
\begin{eqnarray*}
I'_p(u)+\sigma I_p^2(u)=0.
\end{eqnarray*}
\end{theorem}
{\it Proof}. Note that $v=e'(u)$. By Theorem \ref{thmWp}, we have
\begin{eqnarray*}
{1\over 2p t^{a+1}}{d\over dt}\mathcal{W}_p(u)
&=&-\int_M \left[\left\|\nabla^2 v+{g\over (b+2)t}\right\|_{\rm HS}^2+{1\over m-n}\left(\nabla\phi\cdot\nabla v-
{m-n\over (b+2)t}\right)^2\right]u^pd\mu\nonumber\\
& &-\int_M Ric_{m, n}(L)(\nabla v, \nabla v)u^pd\mu-(p-1)\int_M \left(L v+{\kappa\over t}\right)^2u^pd\mu.   \label{WWp-5}
\end{eqnarray*} 
and by Theorem \ref{thm5}, it holds

\begin{eqnarray*}
{\|u\|_p^p\over 2\sigma N_p(u)}{d^2N_p\over dt^2}
&=&- \left(  p-1+{1\over m} \right) \int_M \left[ L v-\int_M L v d\gamma \right]^2u^pd\mu-\int_M Ric_{m, n}(L)(\nabla v, \nabla v) u^pd\mu\nonumber\\
& &\ - \int_M \left[ {m-n\over mn}\left(L v+{m\over m-n}\nabla\phi\cdot\nabla v\right)^2 + \left\|\nabla^2 v-{\Delta v\over n} g\right\|^2_{\rm HS} \right]u^pd\mu, \label{WWp-6}
\end{eqnarray*}
where $d\gamma={u^p d\mu\over \int_M u^p d\mu}={vud\mu\over \int_M vud\mu}$.
%

Let $I_1={\|u\|_p^p\over 2\sigma N_p(u)}{d^2 N_p(u)\over dt^2}$, $I_2={1\over 2p t^{a+1}}{d\over dt}\mathcal{W}_p(u)$ and $I_3=I_1-I_2$. Then $I_3=I_{31}+I_{32}+I_{33}$, where 
\begin{eqnarray*}
I_{31}&=&\int_M \left[ \left\|\nabla^2 v+{g\over (b+2)t}\right\|_{\rm HS}^2-\left\|\nabla^2 v-{\Delta v\over n} g\right\|^2_{\rm HS}\right] u^pd\mu\\
&=&\int_M \left[{2\over (b+2)t}\Delta v+{1\over n}|\Delta v|^2+{n\over (b+2)^2 t^2}\right] u^pd\mu\\
&=&\int_M \left[ {2\over (b+2)t}\Delta v+{1\over n}(Lv+\nabla \phi\cdot \nabla v)^2+{n\over (b+2)^2 t^2}\right] u^pd\mu\\
&=&\int_M \left[ {2\over (b+2)t}\Delta v+{n\over (b+2)^2 t^2}\right] u^pd\mu\\
& &+\int_M \left[ {1\over m}|Lv|^2-{|\nabla \phi\cdot \nabla v|^2 \over m-n}+{m-n\over mn}\left(Lv+{m\over m-n}\nabla\phi\cdot\nabla v\right)^2\right] u^pd\mu,
\end{eqnarray*}
\begin{eqnarray*}
I_{32}&=&\int_M (p-1) \left(L v+{\kappa\over t}\right)^2 u^pd\mu-\left(  p-1+{1\over m} \right) \int_M \left|L v-\int_M L v d\gamma \right|^2 u^pd\mu\\
&=&\int_M \left[-{1\over m}|Lv|^2 +{2a\over t}Lv +{a^2\over (p-1)t^2}\right]u^pd\mu+\left(p-1+{1\over m}\right)\left(\int_M  Lv u^p d\mu\right)^2\left(\int_M u^p d\mu\right)^{-1},
\end{eqnarray*}
and 
\begin{eqnarray*}
I_{33}&=&\int_M \left[ {1\over m-n}\left(\nabla\phi\cdot\nabla v-
{m-n\over (b+2)t}\right)^2-   {m-n\over mn}\left(L v+{m\over m-n}\nabla\phi\cdot\nabla v\right)^2   \right]  u^pd\mu.
\end{eqnarray*}
Substituting $\kappa={m\over  m(p-1)+2}$,  $a=(p-1)\kappa$ and $b=m(p-1)$ into the above calculation, we can derive  
\begin{eqnarray*}
I_3&=&{m(b+1)\over (b+2)^2t^2}\int_M u^pd\mu+{2(b+1)\over (b+2)t}\int_M Lv u^pd\mu +{b+1\over m}\left(\int_M  Lv u^p d\mu\right)^2\left(\int_M u^p d\mu\right)^{-1}\\
&=& {1+m(p-1)\over m}\left[\int_M Lv u^p d\mu+{m\over (2+m(p-1))t}\int_M u^pd\mu \right]^2\left[\int_M u^p d\mu\right]^{-1}\nonumber\\
&=& {1+m(p-1)\over m}\left[\int_M  Lv d\gamma+{m\over (2+m(p-1))t} \right]^2\left[\int_M u^p d\mu\right]\nonumber.
\end{eqnarray*}
Similarly to $(\ref{LvIp})$,  $\int_M Lv d\gamma=-I_p(u)$.  This completes the proof of Theorem \ref{NIW3}. \hfill $\square$

\begin{remark}  In particular, taking limit $p\rightarrow 1$ in Theorem \ref{NIW3}, we recover the NIW formula $(\ref{NIW-1})$ in Theorem \ref{NIW1} for the Shannon entropy power $N(u)$, the Fisher information $I(u)$ and the $W$-entropy associated with 
 the heat equation $\partial_t u=L u$ on Riemannian manifold, i.e., 
\begin{eqnarray*}
 {d^2 N(u)\over dt^2}={2\over m} N(u) \left[  {2\over m}\left|I(u)-{m\over 2t} \right|^2+{1\over  t}{d\over dt}\mathcal{W}(u)\right].
 \end{eqnarray*}
 Under the condition $CD(0, m)$, we have
  \begin{eqnarray*}
 {d\mathcal{W}(u)\over dt}\leq  - {2t \over m}\left|I(u)-{m\over 2t} \right|^2. 
 \end{eqnarray*}
 In particular, under the condition $Ric_{m, n}(L)\geq 0$, the equality in the above inequality holds if and only if ${d^2 N\over dt^2}=0$, which is equivalent to say  that $(M, g)$ is quasi-Ricci flat, i.e., $Ric_{m, n}(L)=0$, and $v=\log u$ satisfies 
 $\nabla^2 v=-{I(u)\over m} g$, 
where $I(u)$ is a solution to the differential equation 
\begin{eqnarray*}
I'(u)+{2\over m} I^2(u)=0.
\end{eqnarray*}

\end{remark}

\medskip

\section{Entropy isoperimetric and Gagliardo-Nirenberg-Sobolev  inequalities}

In \cite{LL20}, we proved an entropy isoperimetric inequality for the  Shannon entropy power on complete Riemanniabn 
 manifolds with non-negative Ricci curvature or the $CD(0, m)$-condition and 
the maximal volume growth condition

\begin{theorem}\label{thm10} Let $(M, g, \mu)$ be a complete Riemannian manifold with $CD(0, m)$-condition and the maximal volume growth condition
\begin{eqnarray}
\mu(B(x, r))\geq C_m r^m, \ \ \ \forall x\in M, r>0.
\end{eqnarray}
Then the entropy isoperimetric inequality holds for the Shannon entropy power: for any smooth probability distribution $fd\mu$ such that 
$I(f)$ and $H(f)$ are well-defined, we have
\begin{eqnarray}
I(f)N(f)\geq \gamma_m:=2\pi e m\kappa_*^{2\over m}.  \label{enisop=1}
\end{eqnarray}
where $H(f)=-\int_M f\log f d\mu$, $I(f)=\int_M  {|\nabla f|^2\over f}d\mu$, and
\begin{eqnarray*}
\kappa_*:=\lim\limits_{r\rightarrow \infty}\inf\limits_{x\in M} {V(B(x, r))\over \omega_m r^m},
\end{eqnarray*}
here $\omega_m$ denotes the volume of the unit ball in $\mathbb{R}^m$ (In the case $m=n$ and $M=\mathbb{R}^n$ with the Lebesgue measure, 
$\kappa_*=1$). Equivalently, the Stam type  logarithmic Sobolev inequality holds: for any smooth  function $f$ such that $\int_M f^2d\mu=1$ and $\int_M |\nabla f|^2 d\mu<\infty$, we have
\begin{eqnarray*}
\int_M f^2\log f^2 d\mu\leq {m\over 2}\log \left({4\over \gamma_m} \int_M |\nabla f|^2d\mu\right). \label{LSI3}
\end{eqnarray*}
\end{theorem}

In particular, the above result holds with $m=n$ on complete Riemannian manifolds with non-negative Ricci curvature and with the maximal volume growth condition 
for the standard Riemannian volume measure, which extends the well known entropy isoperimetric inequality for the Shannon entropy power on 
Eucildean space $\mathbb{R}^n$ with the optimal constant $\gamma_n=2\pi e n$. 

The proof of Theorem \ref{thm9} (even in the case $M=n$, $L=\Delta$ and $\mu=v$) uses 
of the concavity of the Shannon entropy power along the heat equation $\partial_t u=Lu$, two sides heat kernel estimates on complete Riemannian manifolds with 
$CD(0, m)$-condition and the maximal volume growth condition. It would be an interesting question whether we can extend Theorem \ref{thm10} to the Renyi entropy 
power on complete Riemnnian manifolds. As far as we know, it seems to be unknown in the literature. We would like to point out that, even thought it would be possible to extend the concavity of the Renyi entropy power to the  nonlinear diffusion equation $\partial_t u=\Delta u^p$ on 
complete Riemannian manifolds with non-negative Ricci curvature  and to the  nonlinear diffusion equation $\partial_t u=Lu^p$  on weighted 
complete Riemannian manifolds with non-negative $m$-dimensional Bakry-Emert Ricci curvature (i.e., $Ric_{m, n}(L)\geq 0$ for general $m>n$), however it is still 
unknown whether we can have nice asymptotic behavior of the fundamental solution of the nonlinear diffusion equation $\partial_t u=\Delta u^p$ or $\partial_t u=Lu^p$ 
on complete non-compact Riemannian manifolds with non-negative Ricci curvature or non-negative $m$-dimensional Bakry-Emert 
Ricci curvature. Thus, it is not available for us to use the same idea as we used for the
 linear heat equation case in \cite{LL20} to prove the entropy isoperimetric inequality for the Renyi entropy power on 
 complete Riemannian manifolds with non-negative Ricci curavture or non-negative $m$-dimensional Bakry-Emery Ricci curvature and 
 maximal volume growth condition. 

On the other hand, for any $n>2$ and $p>{n\over n+2}$, Savar\'e and Toscani \cite{ST} proved that for  any smooth, strictly positive and rapidly decreasing probability density $f$  
on $\mathbb{R}^n$, the following entropy isoperimetric inequality holds
\begin{eqnarray}
N_p(f)I_p(f)\geq \gamma_{n, p}   \label{EIsop}
\end{eqnarray}
with the optimal constant 
$$\gamma_{n, p}=N_p(M_p)I_p(M_p),$$ where $M_p$ is the so-called Barenblatt 
solution to the nonlinear diffusion equation $\partial_t u=\Delta u^p$ on $\mathbb{R}^n$. More precisely, 
\begin{eqnarray*}
M_p(x, t)={1\over t^{n/\mu}}\bar{M}_p\left({x\over t^{1/\nu}}\right),
\end{eqnarray*}
where $\nu=2+n(p-1)$, and 
\begin{eqnarray*}
\bar{M}_p(x)=\left(C-{1\over 2\mu} {p-1\over p}|x|^2\right)_{+}^{1\over p-1};
\end{eqnarray*}
here $(s)_{+}=\max\limits\{s, 0\}$, and the constant $C$ can be chosen to fix the mass of the source-type  Barenblatt 
solution equal to one. For explicit formulae of the best constants $\gamma_{n, p}$ for $p\in (n/n+2, 1)$ and $p>1$, see \cite{ST}.

According to Toscani \cite{Tos2},  the entropy isoperimetric inequality $(\ref{EIsop})$ on $\mathbb{R}^n$
 is equivalent to the  Gagliardo-Nirenberg-Sobolev inequality on $\mathbb{R}^n$  (see $(\ref{GNS1})$ below).  When $p=1-{1\over n}$, Savar\'e and 
Toscani \cite{ST} proved the equivalence between the entropy isoperimetric inequality $(\ref{EIsop})$  and the Sobolev inequality on  $\mathbb{R}^n$ with 
the best Soblev constant: for $n>2$, 
\begin{eqnarray*}
\int_{\mathbb{R}^n} 
|\nabla g|^2dx\geq \left({n-2\over 2n-2}\right)^2 \gamma_{n, 1-{1\over n}} \left(\int_{\mathbb{R}^n} g^{2n\over n-2}dx\right)^{n-2\over n}. \label{SobEu}
\end{eqnarray*}
For the best constants and the extremal functions in the 
Gagliardo-Nirenberg-Sobolev inequality and the Sobolev inequality on $\mathbb{R}^n$, see \cite{BGL, DD, DT, ST} and references therein.

The purpose of this section is to prove the entropy isoperimetric inequality for the Renyi entropy power 
and the Gagliardo-Nirenberg-Sobolev inequality on complete Riemannian manifolds with non-negative Ricci curavture or non-negative $m$-dimensional Bakry-Emery Ricci curvature and 
 maximal volume growth condition. To do so, we need  a new approach.

The idea of Toscani \cite{Tos2} and Savar\'e-Toscani \cite{ST} can go straightforward to (weighted) complete Riemannian manifolds. In fact,   for $p>{m\over m+2}$, 
the entropy isoperimetric inequality 
\begin{eqnarray}
N_p(f)I_p(f)\geq \gamma_{m, p},   \label{EIIm}
\end{eqnarray}
holds on an $n$-dimensional weighted complete Riemannian manifold $(M, g, \mu)$ if and only if the following 
functional inequality: for any smooth and positive probability density $f$ on $M$, 
\begin{eqnarray}
\int_M {|\nabla f^p|^2\over f}d\mu\geq \gamma_{m, p}\left(\int_M f^p d\mu\right)^{2+{2\over 2m(p-1)}}.  \label{finq1}
\end{eqnarray}
Let $f= g^{2\over 2p-1}$ with smooth $g>0$ and $\int_M g^{2\over 2p-1}d\mu=1$. 
Then $(\ref{finq1})$ is equivalent to 
\begin{eqnarray*}
\int_M |\nabla g|^2d\mu \geq \gamma_{m, p} \left({2p-1\over 2p}\right)^2 \left(\int_M g^{2p\over 2p-1}d\mu\right)^{2+{2\over m(p-1)}}. \label{finq2}
\end{eqnarray*}
Equivalently, for smooth $g>0$ with  $\int_M g^{2\over 2p-1}d\mu=1$, 
\begin{eqnarray}
\|\nabla g\|_2^2 \geq \gamma_{m, p} \left({2p-1\over 2p}\right)^2 \|g\|_{2p\over 2p-1}^{{2p\over 2p-1}\cdot {2+2m(p-1)\over m(p-1)}}.
 \label{finq3}
\end{eqnarray}
By a simply scaling argument,    $(\ref{finq3})$ is equivalent to:  for any smooth $g>0$,
\begin{eqnarray}
\|\nabla g\|_2^2\|g\|_{2\over 2p-1}^
{{2p\over 2p-1}\cdot {2+2m(p-1)\over m(p-1)}-2}\geq \gamma_{m, p}  \left({2p-1\over 2p}\right)^2 
 \|g\|_{2p\over 2p-1}^{{2p\over 2p-1}\cdot {2+2m(p-1)\over m(p-1)}}. \label{finq4}
\end{eqnarray}

The functional inequality $(\ref{finq4})$ is indeed equivalent to the Gagliardo-Nirenberg-Sobolev inequality 
\begin{eqnarray}
\|g\|_{q+1}^2\leq A\|\nabla g\|_2^{\theta} \|g\|_{2q}^{1-\theta},\label{GNS1}
\end{eqnarray}
where $q\in (0, 1)$ and $\theta={m(1-q)\over (1+q)(m-(m-2)q)}\in (0, 1)$ satisfy 
\begin{eqnarray}
{1\over q+1}={\theta\over 2^*}+{1-\theta\over 2q}, \label{theta-GNS}
\end{eqnarray}
in which for $m>2$, 
\begin{eqnarray*}
2^*={2m\over m-2}.\label{2*}
\end{eqnarray*}
Indeed, taking $q={1\over 2p-1}$, we can check that $(\ref{finq4})$  is equivalent to $(\ref{GNS1})$ with 
\begin{eqnarray}
A=\left({2p\over 2p-1}\right)^{\theta}\gamma_{m, p}^{-{\theta\over 2}}. \label{A-GNS}
\end{eqnarray}
That is to say, the entropy isoperimetric inequality $(\ref{EIIm})$ with the optimal constant $\gamma_{m, p}$ for $p>{m\over m+2}$ is equivalent to the  optimal 
Gagliardo-Nirenberg-Sobolev inequality 
$(\ref{GNS1})$ with $q={1\over 2p-1}$, $\theta$ in  $(\ref{theta-GNS})$ and with the optimal constant $A$ in $(\ref{A-GNS})$.

In particular, when $m>2$, for $p=1-{1\over m}$, we have ${2+2m(p-1)\over m(p-1)}=0$. The  Gagliardo-Nirenberg-Sobolev inequality $(\ref{finq4})$ becomes the 
Sobolev inequality
%
%
\begin{eqnarray*}
\int_M |\nabla g|^2d\mu\geq \left({m-2\over 2m-2}\right)^2 \gamma_{m, 1-{1\over m}} \left(\int_M g^{2m\over m-2}d\mu\right)^{m-2\over m}.
\end{eqnarray*}
That is to say, the entropy isoperimetric inequality $(\ref{EIIm})$ with the optimal constant $\gamma_{m, p}$ for  $p=1-{1\over m}$ for $m>2$ is equivalent to the  Sobolev inequality with the optimal Sobolev constant:

\begin{eqnarray}
\|g\|_{2m\over m-2}^2 \leq  {1\over \gamma_{m, 1-{1\over m}} } \left({2m-2\over m-2}\right)^2\|\nabla g\|_2^2. \label{Sobm2}
\end{eqnarray}
Moreover, it is well known that, by Proposition 6.2.3 in \cite{BGL} and references therein, the Sobolev inequality $(\ref{Sobm2})$ is equivalent to the Nash inequality 
\begin{eqnarray*}
\|f\|_2^{1+{2/m}} \leq C_m\|\nabla f\|_2 \|f\|_1^{2/m}.
\end{eqnarray*}

Recall the following result which was proved by the second named author in \cite{Li09}. 

\begin{theorem}\label{thm9} Let $(M, g, \mu)$ be a complete Riemannian manifold with $CD(0, m)$-condition  (i.e., $Ric_{m, n}(L)\geq 0$) and the maximal volume growth condition
\begin{eqnarray} 
\mu(B(x, r))\geq C_m r^m, \ \ \ \forall x\in M, r>0. \label{maxV}
\end{eqnarray}Then the Sobolev inequality holds for the measure $\mu$: there exists a positive constant $C_{\rm Sob}$ depending on $m$ and 
$C_m$ such that 
\begin{eqnarray}
\|f\|_{2m\over m-2}^2 \leq C_{\rm Sob} \|\nabla f\|_2^2, \ \ \ \forall f\in C_0^\infty(M). \label{L2Sob2}
\end{eqnarray}
\end{theorem}

The above result was proved in \cite{Li09} using the upper bound heat kernel estimate of the Witten Laplacian obtained by the second named author in \cite{Li05} on 
complete Riemnnian manifolds and the well known criterion due to Carlen-Kusuoka-Stroock \cite{CKS} on the equivalence between 
the validity of the Sobolev inequality and the ultracontractivity of the heat semigroup. More precisely, by Theorem 5.4 in \cite{Li05},  
under the $CD(0, m)$-condition, for any $\varepsilon>0$, $x, y\in M$ and $t>0$, there is a constant $C_1(\varepsilon)>0$ such that 
\begin{eqnarray}
p_t(x, y)\leq {C_1(\varepsilon)\over\sqrt{ \mu(B(x, \sqrt{t})) \mu(B(y, \sqrt{t}))}} e^{-{d^2(x, y)\over (4+\varepsilon)t}}. \label{heatupperbound}
\end{eqnarray}
Now the maximal volume growth condition $(\ref{maxV})$ yields
\begin{eqnarray}
p_t(x, y)\leq {C_2(\varepsilon)\over t^{m/2}} e^{-{d^2(x, y)\over (4+\varepsilon)t}}. 
\end{eqnarray}
Hence $P_t=e^{tL}$ is ultracontractive in the sense that
\begin{eqnarray}
\sup\limits_{x, y\in M}p_t(x, y)\leq {C_2(\varepsilon)\over t^{m/2}}. \label{ultra}
\end{eqnarray}
By  the well known criterion due to Carlen-Kusuoka-Stroock \cite{CKS}, $(\ref{ultra})$ is equivalent to say that the Sobolev inequality $(\ref{L2Sob2})$ holds. 

We can use the entropy isoperimetric inequality $(\ref{enisop=1})$ to give an alternative proof of Theorem \ref{thm9}. Indeed, by 
Theorem \ref{thm10},  under the condition of Theorem \ref{thm9}, the entropy isoperimetric inequality $(\ref{enisop=1})$ holds. Equivalently, 
the entropy-energy inequality (in the sense of \cite{BGL}) holds: for all $f\in C_0^\infty(M)$ with $\int_M f^2d\mu=1$,  
\begin{eqnarray}
{\rm Ent}(f^2)\leq {m\over 2} \log \left({4\over \gamma_m}\int_M |\nabla f|^2d\mu\right). \label{EE}
\end{eqnarray}
By Proposition 6.2.3 in \cite{BGL},  the entropy-energy inequality $(\ref{EE})$ implies  the Nash inequality
\begin{eqnarray}
\|f\|_2^{1+{2/m}} \leq {2\over \sqrt{\gamma_m}} \|\nabla f\|_2\|f\|_1^{2/m}, \ \ \ \forall f\in C_0^\infty(M), \label{Nash0}
\end{eqnarray}
which is equivalent to the Sobolev inequality $(\ref{L2Sob2})$ for some constant $C_{\rm Sob}$ depending only on $m$ and $\gamma_m$. Indeed, all of 
$(\ref{EE})$, $(\ref{Nash0})$ and $(\ref{L2Sob2})$ are equivalent to 
the $L^1$ to $L^\infty$-ultracontractivity of the heat semigroup $P_t=e^{tL}$, i.e., 
\begin{eqnarray}
\|P_tf\|_{\infty}\leq {C\over t^{m/4}}\|f\|_2, \ \ \ \ \forall t>0, \label{ulta0}
\end{eqnarray}
or the $L^1$ to $L^\infty$-ultracontractivity  of  the heat semigroup $P_t=e^{tL}$, i.e., 
\begin{eqnarray}
\|P_tf\|_{1, \infty}\leq {C\over t^{m/2}}\|f\|_1, \ \ \ \ \forall\ t>0. \label{ulta1}
\end{eqnarray}

Moreover, by Proposition 6.10.2 in \cite{BGL}, the Sobolev inequality $(\ref{L2Sob2})$ is equivalent to the following Gagliardo-Nirenberg-Sobolev inequality 
\begin{eqnarray}
\|f\|_{q}\leq C^{\theta\over 2}_{\rm Sob}\|\nabla f\|_2^{\theta} \|f\|_{s}^{1-\theta}, \ \ \ \ \forall f\in C_0^\infty(M)\label{GNS2}
\end{eqnarray}
for any $1\leq s\leq  q \leq 2^*= {2m\over m-2}$, where $q$, $s$ and $\theta$ satisfy 
\begin{eqnarray}
{1\over q}={\theta\over 2^*}+{1-\theta\over s}. \label{theta-GNS2}
\end{eqnarray}
In particular, when $q=2$, $s=1$ and $\theta={m\over m+1}$, $(\ref{GNS2})$ becomes the Nash inequality
\begin{eqnarray}
\|f\|_{2}^{1+2/m}\leq \sqrt{C_{\rm Sob}}\|\nabla f\|_2 \|f\|_{1}^{2/m}, \ \ \ \ \forall f\in C_0^\infty(M)\label{GNS2}
\end{eqnarray}

In view of this argument and Theorem \ref{thm9}, we have the following  result which gives a reasonable condition for the validity of the entropy isoperimetric 
inequalities and the  Gagliardo-Nirenberg-Sobolev inequality on complete Riemannian manifolds. As far as we know, it is new in the literature even in the case $m=n$, 
$\phi=0$, $Ric_{m, n}(L)=Ric$, and $\mu=v$. 

\begin{theorem}\label{thm11} Let $(M, g, \mu)$ be a complete Riemannian manifold with $CD(0, m)$-condition (i.e.,  $Ric_{m, n}(L)\geq 0$) and the maximal volume growth condition $(\ref{maxV})$, i.e., 
\begin{eqnarray*}
\mu(B(x, r))\geq C_m r^m, \ \ \ \forall x\in M, r>0.
\end{eqnarray*}
Then\\
$(i)$ the entropy isoperimetric inequality $(\ref{enisop=1})$ holds, i.e., 
\begin{eqnarray*}
I(f)N(f)\geq \gamma_m:=2\pi e m\kappa_*^{2\over m}.  
\end{eqnarray*}
Moreover, there exists a positive constant $C_{\rm Sob}$ depending on $m$ and $\gamma_m$ such that \\
$(ii)$ the Sobolev inequality $(\ref{L2Sob2})$ holds, i.e.,  for any $f\in C_0^\infty(M)$, 
\begin{eqnarray*}
\|f\|_{2m\over m-2}^2 \leq C_{\rm Sob} \|\nabla f\|^2.
\end{eqnarray*}
$(iii)$ the  Gagliardo-Nirenberg-Sobolev inequality holds: for any $f\in C_0^\infty(M)$, 
\begin{eqnarray*}
\|f\|_{q}\leq C^{\theta\over 2}_{\rm Sob}\|\nabla f\|_2^{\theta} \|f\|_{s}^{1-\theta}, \ \ \ \ \forall f\in C_0^\infty(M)\label{GNS2}
\end{eqnarray*}
for any $1\leq s\leq  q \leq 2^*= {2m\over m-2}$, where $q$, $s$ and $\theta\in [0, 1]$ satisfy 
\begin{eqnarray*}
{1\over q}={\theta\over 2^*}+{1-\theta\over s}. \label{theta-GNS2}
\end{eqnarray*}
$(iv)$ the Nash inequality holds,  i.e., for any $f\in C_0^\infty(M)$, 
\begin{eqnarray*}
\|f\|_2^{1+{2/m}} \leq \sqrt{C_{\rm Sob}} \|\nabla f\|_2\|f\|_1^{2/m}.
\end{eqnarray*}
$(vi)$ the entropy isoperimetric inequality holds for the Renyi entropy: for any smooth probability distribution $fd\mu$ such that 
$I_p(f)$ and $H_p(f)$ are well-defined, we have
\begin{eqnarray*}
N_p(f)I_p(f)\geq \gamma_{m, p}.
\end{eqnarray*}
where $\gamma_{m, p}$ is a constant depending on $m$, $p$ and $C_{\rm Sob}$.\\
$(vii)$ the ultracontractivity of the heat semigroup holds: there exists a constant $C>0$ depending on $m$ and $C_{\rm Sob}$ such that 
\begin{eqnarray*}
p_t(x, y)\leq {C\over t^{m/2}}, \ \ \ \forall x, y\in M, t>0.
\end{eqnarray*}
\end{theorem}

\begin{remark} In \cite{Deman1, Deman2, Deman3}, Demange proved some Sobolev inequalities  on complete Riemannian manifolds with 
Ricci curvature bounded below by zero, positive or negative constant. 
More precisely, let $M$ be a $d$-dimensional complete Riemannian manifold with Ricci curvature bounded below by a negative constant $-\rho$, where $\rho>0$, and
 assume that there is a constant $\beta\in \mathbb{R}$ such that  the potential function $T: M\rightarrow \mathbb{R}$ satisfies  
\begin{eqnarray*}
\nabla^2 T\geq \left({\rho\over d-1}T+\beta\right)g. 
\end{eqnarray*}
Let $v>0$ be such that $(d-1)v^{-1/d}=T$ and suppose that the following technical condition holds: there exists $m\geq d$ 
such that 
\begin{eqnarray*}
\int (T^{-m}|\nabla T|^4 +T^{-d})d\mu<\infty.
\end{eqnarray*}
 Then for any smooth positive function $f$ with $\int_M fd\mu=\int vd\mu$ (where $\mu$ denotes the volume measure) and $f=v$ on the complementary of
 a compact set, the following Sobolev inequality holds
\begin{eqnarray*}
2\beta E(f)+{\rho\over d-1}J(f)\leq I(f),
\end{eqnarray*}
where $I(f)=\int_M f|\nabla(-(d-1)f^{-1/d}+T)|^2d\mu$, $J(f)=\int_M \left[\int_{v(x)}^{f(x)} (-(d-1)s^{-1/d}+T(x))^2ds \right]d\mu(x)$ and 
$E(f)=\int_M \left[\int_{v(x)}^{f(x)} (-(d-1)s^{-1/d}+T(x))ds \right]d\mu(x)$.  See \cite{Deman1}. 
On the other hand, assuming that $M$ is a $d$-dimensional complete Riemannian manifold with nonnegative Ricci curvature and equipped with a potential function $T=(d-1)v^{-1/d}$ such that for 
a strictly positive constant $\rho>0$, $\nabla^2 T\geq \rho g$. Then the following Sobolev inequality holds
\begin{eqnarray*}
\int_M \left(d{v^{1-1/d}-f^{1-1/d}
\over (d-1)^2}+v^{-1/d}{f-v\over d-1}\right)d\mu\leq {1\over 2\rho}\int_M f|\nabla(f^{-1/d}-v^{-1/d})|^2d\mu.
\end{eqnarray*}
See \cite{Deman3}. The case of compact Riemannian manifolds with positive Ricci curvature lower bound has been also studied. See \cite{Deman2}. 
\end{remark}
\begin{remark} {\ As we have pointed out in Introduction and in the beginning of Section $4$,  we can use a similar argument as used in  the first proof of 
the entropy power concavity inequalities  (EPCI) in Theorem \ref{thm1}, Theorem \ref{thm2} and Theorem \ref{thm3}  to prove the EPCI for the Renyi entropy power 
associated with  nonlinear diffusion equation on RCD$(K, N)$ metric measure 
spaces and on compact metric measure spaces equipped with a $(K, N)$-super Ricci flow in the sense of \cite{Sturm18}. By the same argument as used in the proof of Theorem \ref{thm11}, we can further extend the entropy isoperimetric inequalities for the 
Shannon entropy power and the Renyi entropy power,  the Gagliardo-Nirenberg-Sobolev inequality, the Sobolev inequality and the Nash inequality 
to RCD$(0, N)$ spaces with the maximal volume growth condition. This will be developed in a forthcoming paper.  

}
\end{remark}

To end this paper, we would like  to mention that Wang et al. \cite{W1, W2} proved the concavity of the Renyi entropy power for the
$p$-Laplacian equation  and the double nonlinear diffusion equations on compact Riemannian manifolds. Our work on 
the entropy power concavity inequalities in Theorem \ref{thm1}, Theorem \ref{thm2} and 
Theorem \ref{thm3} has been  done in our 2017 preprint \cite{LL17} and is independent of \cite{W1, W2}.  
\medskip

\noindent{\bf Acknowledgement}. The second author is very grateful to Prof. N. Mok for his suggestion to study the entropy power on manifolds 
and Prof. L. Ni for his suggestion to study nonlinear diffusion equations on manifolds.  We would like to thank Prof. Guang-Yue Han and Prof. Feng-Yu Wang for helpful discussions in the  earlier stage 
of  this work. Finally, we would like to thank  Dr. Yu-Zhao Wang for his interest, nice  
discussions  and for his careful reading on the earlier versions of the manuscript of this work. 
%
%


\begin{flushleft}
\medskip\noindent

Songzi Li, School of Mathematics,  Renmin University of China,  Beijing, 100872, China
Email: sli@ruc.edu.cn
\medskip

Xiang-Dong Li, Academy of Mathematics and Systems Science, Chinese
Academy of Sciences, No. 55, Zhongguancun East Road, Beijing, 100190,  China\\
E-mail: xdli@amt.ac.cn

and

School of Mathematical Sciences, University of Chinese Academy of Sciences, Beijing, 100049, China
\end{flushleft}


\begin{thebibliography}{99}
\bibitem{BE} D. Bakry, M. Emery, Diffusion hypercontractives, S\'em. Prob. XIX, Lect. Notes in Maths. 1123 (1985), 177-206.
\bibitem{BGL} D. Bakry, Y. Gentil, M. Ledoux, Analyse and Geometry for Markov Diffusion Operators, Springer, 2014. 
\bibitem{Bla} N.M. Blachman, The convolution inequality for entropy powers. IEEE Trans. Inform. Theory 2, 267-271, (1965).
\bibitem{CKS} E. Carlen, S. Kusuoka, D. Stroock, Upper bounds for symmetric Markov transition functions. Ann. Inst. H. Poincar\'e, Probab. Stat. 23 (1987), 245-287.
\bibitem{Bess} A. L. Besse, Einstein manifolds, Ergebnisse der Mathematik (3) 10, Springer, Berlin, 1987.
\bibitem{Cost} M.H.M. Costa, A new entropy power inequality, IEEE Trans. Inf. Theory, IT- 31, (6) 751-760, (1985).
\bibitem{Dav} E.B. Davies, Heat Kernel and Spectral Theory, Cambridge University Press, Cambridge, New York, 1990.
\bibitem{DD} M. Del Pino, J. Dolbeault, Best constants for Gagliardo-Nirenberg inequalities and applications to nonlinear diffusions,
J. Math. Pures Appl. (9) 81 (2002) 847-875.
\bibitem{Deman1} J. Demange,  Porous media equation and Sobolev inequalities under negative curvature, Bull. Sci. Math. 129 (2005) 804-830.
\bibitem{Deman2} J. Demange,  Improved Gagliardo–Nirenberg–Sobolev inequalities on manifolds with positive curvature, J. Funct. Anan. 254 (2008) 593-611.
 \bibitem{Deman3} J. Demange, From porous media equation to generalized Sobolev inequalities on a Riemannian manifold, preprint. 2004.
\bibitem{DT} J. Dolbeault, G. Toscani, Nonlinear diffusions: Extremal properties of Barenblatt profiles, best matching and delays, Nonlinear Analysis
38 (2016), 31-43.
\bibitem{Dem1} A. Dembo, A simple proof of the concavity of the entropy power with respect to the variance of additive normal noise, IEEE Trans. Inform. Theory 35, 887-888, (1989).
\bibitem{Dem2}  A. Dembo, T.M. Cover, and J.A. Thomas, Information theoretic inequalities, IEEE Trans. Inf. Theory, 37, (6), 1501-1518, (1991).
\bibitem{HL} G. Huang, H. Li,  Gradient estimates and entropy formulae of porous medium and fast diffusion equations for the Witten Laplacian, Pacific J. Math.  268 (2014), No. 1, 47-78.
\bibitem{KL} K. Kuwada, X.-D. Li,  Monotonicity and rigidity of the W-entropy on RCD(0,N) spaces, Manuscripta Mathematica 2020. https://doi.org/10.1007/s00229-019-01177-y
\bibitem{Li05} X.-D. Li,  Liouville theorems for symmetric diffusion operators on complete Riemannian manifolds, J. Math. Pures Appl. 84 (2005), 1295-1361.
\bibitem{Li09} X.-D. Li, Sobolev inequalities on forms and $L^{p,q}$-cohomology on complete Riemannian manifolds, J. Geom. Anal. 20 (2010) 354-387.
\bibitem{Li12}X.-D. Li, Perelman's entropy formula for the Witten Laplacian
on Riemannian manifolds via Bakry-Emery Ricci curvature, Math. Ann. 353 (2012), 403-437.
\bibitem{Li16} X.-D. Li, Hamilton's Harnack inequality and the $W$-entropy formula on complete Riemannian manifolds, 
Stoch. Processes  Appl. 126 (2016) 1264-1283.
\bibitem{LL15}S. Li, X.-D. Li,  $W$-entropy formula for the Witten Laplacian on manifolds with time dependent metrics and potentials,  Pacific J. Math. Vol. 278 (2015), No. 1, 173-199.
\bibitem{LL16}  S. Li, X.-D. Li,  $W$-entropy formulas and Langevin deformation of flows on Wasserstein space over Riemannian manifolds.
ArXiv:1604.02596, 2016.
\bibitem{LL-AJM} S. Li, X.-D. Li, On Harnack inequalities for Witten Laplacian on Riemannian manifolds with super Ricci flows. Asian J
Math, (22) 2018, 577-598.
\bibitem{LL17} S. Li, X.-D. Li, Entropy differential inequality and entropy power inequality on Riemannian manifolds, preprint, 2017.
\bibitem{LL18a} S. Li, X.-D. Li, Hamilton differential Harnack inequality and $W$-entropy for Witten Laplacian on Riemannian manifolds, J. Funct. Anal. 274 (2018) 3263-3290.
\bibitem{LL18b} S. Li, X.-D. Li, $W$-entropy formulas on super Ricci flows and Langevin deformation on Wasserstein space over Riemannian
manifolds, Sci China Math, (61) 2018, 1385-1406.
\bibitem{LL19a} S. Li, X.-D. Li, On the Li-Yau-Hamilton Harnack inequalities on Ricci flow and super Ricci flows (in Chinese), Sci Sin Math, (49) 2019, no.~11,  1613-1632.
\bibitem{LL19b} S. Li, X.-D. Li, $W$-entropy, super Perelman Ricci flows and $(K, m)$-Ricci solitons,  J Geom Anal, (2013) 2019, 1-32.
\bibitem{LL20} S. Li, X.-D. Li, On the Shannon entropy power on Riemannian manifolds and Ricci flow, arxiv2001.00410v1, 2020.
\bibitem{LY} P. Li, S.-T. Yau, On the parabolic kernel of the Schr\"odinger operator, Acta. Math. 156 (1986), 153-201.
\bibitem{LX} J. Li, X. Xu, Differential Harnack inequalities on Riemannian manifolds I: linear heat equation, Adv. Math. 226 (5) (2011), 4456-4491.
\bibitem{Lot1}  J. Lott, Optimal transport and Perelman’s reduced volume. Calc Var Partial Differential Equations, 36 (2009), 49-84 
\bibitem{LV} J. Lott, C. Villani,  Ricci curvature for metric-measure spaces via optimal transport. Ann of Math. 169  (2) 2009, 903-991
\bibitem{LNVV} P. Lu, L. Ni, J.-L. Vazquez, C. Villani, Local Aronson-B\'enilan estimates and entropy formulae for porous medium and fast diffusion equations on manifolds, J. Math. Pures Appl., Vol 91 (2009), no. 1, 1-19.
\bibitem{N1} L. Ni, The entropy formula for linear equation, J. Geom. Anal. 14 (1), 87-100, (2004).
\bibitem{N2} L. Ni, Addenda to ``The entropy formula for linear equation'', J. Geom. Anal. 14 (2), 329-334, (2004).
\bibitem{Ot} F. Otto, The geometry of dissipative evolution equations: the porous medium equation, Commun. in Parial Differential Equations 26 (1 and 2), 101-174 (2001).
\bibitem{P1} G. Perelman, The entropy formula for the Ricci flow and its geometric applications, http://arXiv.org/abs/maths0211159.
\bibitem{ST}G. Savar\'e, G. Toscani,  The concavity of R\'enyi entropy power,  IEEE TRANSACTIONS ON INFORMATION THEORY,  2014.
\bibitem{Sh48}  C. E. Shannon,  A mathematical theory of communication, Bell Syst. Tech. J., vol. 27, pp. 623-656, Oct. 1948.
\bibitem{Stam}  A. J. Stam, Some inequalities satisfied by the quantities of information of Fisher and Shannon, Information and Control, vol. 2, pp. 101-112, Jun. 1959.
\bibitem{Sturm18} K.-Th. Sturm, Super-Ricci flows for metric measure spaces, 
J. Funct. Anal. 275  (2018), no.~12, 3504--3569.
\bibitem{Tos2} G. Toscani, The information-theoretic meaning of Gagliardo-Nirenberg type inequalities. Atti Accad. Naz. Lincei Rend. Lincei Mat. Appl. 30 (2019), 
no. 2, 237-253
\bibitem{Vaz} J. L. Vazquez, The Porous Medium Equation: Mathematical Theory, Oxford, UK: Oxford University Press, 2007.
\bibitem{V0} C. Villani, A Short Proof of the Concavity of Entropy Power,  IEEE TRANSACTIONS ON INFORMATION THEORY, VOL. 46, NO. 4,1695-1696, July 2000.
\bibitem{V1} C. Villani,  Topics in Mass Transportation, Grad. Stud. Math., Amer. Math. Soc., Providence,
RI, 2003.
\bibitem{V2}C. Villani, Optimal Transport, Old and New, Springer, 2008.
\bibitem{W1} Y.-Z. Wang, X. Zhang, The concavity of $p$-entropy power and applications in functional inequalities, Nonlinear Anal. 179 (2019), 1-14.
\bibitem {W2} Y.-Z. Wang, Y.-M. Wang, The concavity of $p$-R\'enyi entropy power for doubly nonlinear diffusion equations and $L^p$-
Gagliardo-Nirenberg-Sobolev inequalities. J. Math. Anal. Appl. 484 (2020), no. 1, 123698.
\end{thebibliography}
\end{document}